\documentclass{amsart}

\usepackage{amssymb,eepic,epic}

\newtheorem{theorem}{Theorem}[section]
\newtheorem{proposition}[theorem]{Proposition}
\newtheorem{lemma}[theorem]{Lemma}
\newtheorem{corollary}[theorem]{Corollary}

\theoremstyle{remark}
\newtheorem{remark}[theorem]{Remark}

\theoremstyle{definition}

\newcommand{\Ind}{\operatorname{Ind}}
\newcommand{\RR}{\operatorname{RR}}
\newcommand{\Td}{\operatorname{Td}}
\newcommand{\Ch}{\operatorname{Ch}}
\newcommand{\Spin}{\operatorname{Spin}}
\newcommand{\Eul}{\operatorname{Eul}}

\newcommand\lie[1]{\mathfrak{#1}} 
\renewcommand\a{\lie{a}}

\newcommand\g{\lie{g}}

\renewcommand{\H}{\lie{h}}
\newcommand{\K}{\lie{k}}
\newcommand{\G}{\lie{g}}
\newcommand{\T}{\lie{t}}

\newcommand{\N}{\mathbb{N}}
\newcommand{\R}{\mathbb{R}}
\newcommand{\C}{\mathbb{C}}
\newcommand{\Z}{\mathbb{Z}}
\newcommand{\Q}{\mathbb{Q}}

\newcommand{\Sig}{\Sigma}
\newcommand{\bib}{\bibitem}
\newcommand{\lra}{\longrightarrow}
\newcommand{\hra}{\hookrightarrow}
\newcommand{\ra}{\rightarrow}
\newcommand{\A}{\mathcal{A}}

\newcommand{\dirac}{/\kern-1.2ex\partial} 

\renewcommand{\d}{{\mbox{d}}}

\newcommand{\f}{\frac}
\newcommand{\p}{\partial}
\renewcommand{\l}{\langle}
\renewcommand{\r}{\rangle}

\newcommand{\ti}{\tilde}

\begin{document}

\title{Symplectic Surgery and the Spin$^c$-Dirac 
Operator}
\author{Eckhard Meinrenken}
\address{Massachusetts Institute of Technology, Department of
Mathematics, Cambridge, Massachusetts 02139-4307}
\email{mein@math.mit.edu}
\thanks{Revised version, November 1995. 
To appear in {\bf Advances in Mathematics}.}
\date{March 1995}

\begin{abstract}
{Let $G$ be a compact connected Lie group, and $(M,\omega)$ 
a compact Hamiltonian $G$-space, with moment map $J:M\ra\G^*$. 
Under the assumption 
that these data are pre-quantizable, one can construct an associated  
$\Spin^c$-Dirac operator $\dirac_{\C}$, whose equivariant index 
yields a virtual 
representation of $G$. We prove a conjecture of Guillemin and Sternberg 
that 
if $0$ is a regular value of $J$, the 
multiplicity $N(0)$ of the trivial representation in the index 
space $\text{ind}(\dirac_{\C})$, 
is equal to the index of the $\Spin^c$-Dirac operator 
for the symplectic quotient 
$M_0=J^{-1}(0)/G$. This generalizes previous results for the 
case that $G=T$ is abelian, i.e. a torus.}  
\end{abstract}
\maketitle
\tableofcontents

\section{Introduction}
Let $(M,\omega)$ be a compact 
symplectic manifold, with integral symplectic form 
$[\omega]\in H^2(M,\Z)$, and 
$L\ra M$ a Hermitian line bundle whose first 
Chern class is $c_1(L)=[\omega]$. It is then possible to choose a Hermitian 
connection $\nabla$ on $L$ with curvature $F(L)=\f{2\pi}{i}\omega$. 
Choose an $\omega$-compatible, 
positive almost complex structure on $M$, and let 
\begin{equation}\Lambda^\bullet\, T^*M\otimes \C=\oplus_{i,j}\,
\Lambda^{i,j}\,T^*M\end{equation}
be the associated bigrading of the bundle of forms. Using a 
Hermitian connection on the canonical bundle, one has the twisted 
$\Spin^c$-Dirac operator
\begin{equation}\dirac_{\C}:\,\A^{0,even}(M,L)\ra \A^{0,odd}(M,L)\end{equation}
where 
\begin{equation}\A^{i,j}(M,L)=C^\infty(M,\Lambda^{i,j}\,T^*M\otimes L)\end{equation}
is the space of $L$-valued forms of type $(i,j)$.  
The Riemann-Roch number of $L\to M$ is defined to be the dimension of the 
virtual 
vector space
\begin{equation}\text{ind}(\dirac_{\C})=[\text{ker}(\dirac_{\C})]-[\text{ker}
(\dirac_{\C}^*)].\end{equation}
By the index theorem of Atiyah-Singer,
\begin{equation}\RR(M,L)=\int_M \Td(M)\,\Ch(L),\end{equation}
where $\Td(M)$ and $\Ch(L)$ are the Todd class of $M$ and 
Chern character of
$L$, respectively. Since $\Ch(L)=e^\omega$, and since any two 
compatible 
almost complex structures are homotopic, $\RR(M,L)$ is a symplectic 
invariant of 
$M$.

If $M$ is in fact K\"ahler, $L$ a holomorphic Hermitian line bundle, and 
$\nabla$ its canonical connection, then $\dirac_{\C}$ is the Dirac 
operator 
$\sqrt{2}(\bar{\partial}+\bar{\partial}^*)$ associated to the twisted 
Dolbeault complex. If $L$ is  moreover sufficiently positive, 
$\RR(M,L)$ is simply the dimension of the space of holomorphic sections. 
This 
happens for instance for coadjoint orbits $\mathcal{O}=G.\lambda$ 
of compact Lie groups $G$, where $\lambda$ is a dominant weight.
Here $\omega$ is the canonical invariant K\"ahler form on $\mathcal{O}$, and 
$L$ the pre-quantum line bundle. By the Borel-Weil-Bott 
theorem, the space of holomorphic sections of $L$ is just the irreducible 
representation space with highest weight $\lambda$.  

Suppose now that a compact Lie group $G$ acts on $M$ in a Hamiltonian 
fashion, with equivariant moment map $J:M\ra \G^*$. That is,  
\begin{equation}\iota(\xi_M)\omega=\d\l J,\xi\r\end{equation}
for all $\xi\in\g$ where $\xi_M$ is the fundamental vector field.  
Assume also that the action lifts 
to an action on $L$, in such a way that the fundamental vector fields 
on $m$ and $L$ are related by Kostant's formula
\begin{equation}\xi_L=\text{Lift}(\xi_M)+2\pi  \l J,\xi\r\, \f{\p}{\p \phi},
\label{lift}\end{equation}
where ``$\text{Lift}$'' is the horizontal lift with respect to the 
connection and 
$\f{\p}{\p \phi}$ the generating vector field for the scalar 
$S^1$-action 
on $L$. By making the above choices $G$-invariant, we obtain 
$G$-representations 
on $ \text{ker}(\dirac_{\C})$ and $\text{ker}(\dirac_{\C}^*)$, hence a 
virtual 
representation on $\text{ind}( \dirac_{\C})$. The 
cha\-rac\-ter $\RR(M,L)=\chi\in R(G)$ of 
this representation is called the equivariant Riemann-Roch number, and 
the equivariant index theorem of Atiyah-Segal-Singer expresses $\chi(g)$
as an integral of characteristic classes over the fixed point 
manifold 
$M^g$. 
Consider the decomposition of $\chi$ into irreducible characters
$\rho_\mu$ for $G$, 
\begin{equation}\chi=\sum_{\mu\in\hat{G}} 
N(\mu)\,\rho_\mu,\,\,N(\mu)\in\Z.\end{equation}
Let $T$ be the maximal torus of $G$, $\G=\T\oplus [\T,\G]$ the 
corresponding 
decomposition of the Lie algebra, and $\Lambda\subset \T$ the integral 
lattice. 
Let $\T^*_+\subset \G^*$ be some choice of a positive Weyl chamber, 
$\lie{R}_+$ the associated system of positive roots,
and label the irreducible representations 
by the set of dominant weights $ \Lambda^*_+=\Lambda\cap \T^*_+$. 
If $0$ is a regular value of $J$, the action of $G$ on $J^{-1}(0)$ is 
locally 
free, and therefore the symplectically reduced space 
$M_{red}=J^{-1}(0)/G$ is a symplectic orbifold. Moreover, 
$L_{red}:=(L|J^{-1}(0))/G$ is a pre-quantum orbifold-line bundle. 
One has an associated  
$\Spin^c$-Dirac operator, with index $\RR(M_{red},L_{red})$ 
given by the orbifold-index theorem of 
Kawasaki \cite{K79}. 
This paper will be concerned with the proof of a conjecture of 
Guillemin and 
Sternberg \cite{GS82a}, 
that the multiplicity $N(0)=:\RR(M,L)^G$ of the trivial 
representation is exactly 
the Riemann-Roch number $\RR(M_{red},L_{red})$.

Since orbifolds will play an essential 
role in this paper, we will allow that $M$ itself is an 
orbifold, and prove the 
following: 

\begin{theorem}\label{GSC} 
Let $G$ be a compact connected Lie group, and 
$(M,\omega)$ a quantizable Hamiltonian 
$G$-orbifold, with moment map $J:M\ra \G^*$. 
If $0$ is a 
regular value of $J$, 
the multiplicity of the trivial representation is given by 
\begin{equation}\RR(M,L)^G=\RR(M_{red},\,L_{red}).
\label{mult0}\end{equation}
If a compact, connected Lie group $H$ acts on $L\ra M$, 
such that the action 
commutes with the action of $G$, this equality  
holds as an equality of virtual characters for $H$.
\end{theorem} 
By the ``shifting-trick'', one has pre-quantum orbifold-bundles 
$L_\mu$ for all reduced spaces $M_\mu=J^{-1}(\mu)/G_\mu$, 
whenever $\mu\in\Lambda^*_+$ 
is a regular value of $J$, and obtains the following  

\begin{corollary}\label{SupportMultiplicity} 
If  $\mu\in\Lambda^*_+$ is a regular value of $J$, 
$N(\mu)=\RR(M_\mu,L_\mu)$. In 
particular, the support of the multiplicity function is contained in the 
moment polytope $\Delta=J(M)\cap \T^*_+$.
\end{corollary}
For the  K\"ahler case, a variant of Theorem \ref{GSC} was proved by 
Guillemin and 
Sternberg in their 1982 paper, \cite{GS82a}. For the case that 
$G$ is abelian, but $M$ not necessarily K\"ahler, Theorem
\ref{GSC} was proved by Guillemin \cite{G93} under some 
additional assumptions and by Vergne \cite{V94} and Meinrenken \cite{M94}
in general.

For $G$ nonabelian, it was shown in \cite{M94} that Theorem \ref{GSC}
is true if one replaces $L$ by a suitable high tensor power, $L^m$.
This was deduced from a simple stationary phase argument, for which 
it is in fact irrelevant that 
$L$ is a pre-quantum line bundle or that $M$ is
symplectic. For $G=SU(2)$ and under an additional condition on the
moment map, Theorem \ref{GSC} was proved in Jeffrey-Kirwan
\cite{JK94}.  We will show in this paper that neither this extra
condition nor taking tensor powers is necessary, and that Theorem
\ref{GSC} holds for any compact group $G$, without additional hypotheses
besides $0$ being a regular value.

Our approach is based on the symplectic cutting technique of 
E. Lerman \cite{L93}, which was already used in \cite{DGMW95} to 
simplify the proof for the abelian case. 
The first step will be to cut $M$ into smaller pieces, and to prove a 
gluing formula which expresses the Riemann-Roch number of $M$ 
in terms of the Riemann-Roch number of the cut spaces. (Indeed, these 
gluing formulas are valid for arbitrary Hermitian vector bundles, not 
only pre-quantum line bundles.) By combining this result with a ``cross 
section theorem'', we prove the excision property of $N(0)$, i.e. that $N(0)$ 
depends only on local data near $J^{-1}(0)$. This allows us to replace 
$M$ by a new quantizable space $M_{cut}$, such that $M_{cut}$ is 
isomorphic 
to $M$ near $J^{-1}(0)$, and such that the moment polytope of 
$M_{cut}$ is 
contained in a small neighborhood of zero. It can be arranged that 
$M_{cut}$ has 
has very few $T$-fixed point manifolds. By 
explicitly writing down the fixed point formula for 
$\RR(M_{cut},L_{cut})$, we can compute the multiplicity of $0$ 
for $M_{cut}$, and find that it is equal to $\RR(M_{red},L_{red})$.

For the group $G=SU(2)$, an alternative and much easier proof of 
Theorem \ref{GSC} is available. We sketch this proof in the appendix, 
which is independent of the rest of the paper.

If $0$ is {\em not} a regular value of the moment map, the reduced 
space $M_0=\Phi^{-1}(0)/G$ is in general not an orbifold, but 
has more serious singularities. The extension of Theorem \ref{GSC}
to this case will be discussed in a seperate paper with Reyer 
Sjamaar.
\\

\noindent{\bf Acknowledgments.} I would like thank E. Lerman for 
explaining to me his unpublished $SU(2)$-version of symplectic cutting, 
and H. Duistermaat, R. Sjamaar, C. Woodward and M. Vergne 
for helpful comments on an earlier version of this paper.     
This work was supported by a Feodor-Lynen fellowship from 
the Humboldt foundation.

\section{Orbifolds}

In this section, we review some background material on orbifolds, 
and describe the orbifold version of the Berline-Vergne localization 
formula.
Let us briefly recall some basic definitions (for details, see 
Satake \cite{S57} or Kawasaki \cite{K78,K79}). 

Let $M$ be a paracompact Hausdorff topological space.
An orbifold chart for $M$ is a triple $(U,V,H)$ consisting of 
an open subset $U$ of $M$, a finite group $H$, 
an open subset $V$ of some manifold, and a homeomorphism 
$U=V/H$. An orbifold structure on $M$ is a collection of 
orbifold charts $\{(U_i,V_i,H_i)\}$, such that the $U_i$ cover 
$M$ and are subject to appropriate compatibility conditions. 
In particular, one assumes that for 
any two orbifold charts $(U_i,V_i,H_i)$ ($i=1,2)$ and all 
$x\in U_1\cap U_2$, there exists an open neighborhood 
$U\subset U_1\cap U_2$ of $x$ and an orbifold chart $(U,V,H)$, 
together with injections 
$\rho_i:\,H\ra H_i$ and 
$\rho_i$-equivariant 
embeddings $\phi_i:\,V\ra V_i$, with the following property: 
\begin{equation}h_i.\phi_i(V)\cap \phi_i(V)\not=\emptyset
\Rightarrow h_i\in\rho_i(H).\end{equation}
Our definition differs slightly from \cite{K78}, 
in that we do not assume that the actions of the $H_i$ are effective. The 
compatibility conditions ensure that on each connected 
component of $M$, the generic stabilizers $H_{i,*}$ for the 
$H_i$-actions are isomorphic. 
Their order therefore defines a locally constant function on 
$M$, which is called the multiplicity function $d_M:\,M\ra\N$.  

For all $x\in M$, there is a compatible orbifold chart $(U_x,V_x,H_x)$ 
such that the preimage of $x$ in $V_x$ is a fixed point. 
The corresponding $H_x$ is uniquely determined up to isomorphism, 
and is called the isotropy group of $x$. 
Given any finite group $H$, the connected components
of the set of all $x$ such that $H_x\sim H$ are smooth manifolds. 
In this way, 
one obtains a decomposition $M=\cup M_i$, such that $M_i$ is a connected 
manifold and consists of points of a fixed isotropy group, $H_i$. 
Moreover, for all 
$i\not= j$, $M_i$ meets the closure of $M_j$ only if $H_j$ is 
an (abstract) subgroup of $H_i$.   

On each connected component of $M$, there is unique open, dense,  
connected stratum (called principal stratum) 
$M_*$ on which the corresponding isotropy group 
$H_*$ is minimal, i.e. $\# H_*=d_M$. 

If $G$ is a compact Lie group and $G\times P\ra P$ a locally free 
action on a manifold $P$, the orbit space $M=P/G$ is an 
orbifold\footnote{In 
fact, any orbifold is of this form; for example one can take 
$P$ to be 
the $O(n)$ frame bundle with 
respect to some Riemannian metric on $M$, and $H=O(n)$.}.
Indeed, by the slice theorem a neighborhood of any orbit $G.p=x$ is 
equivariantly diffeomorphic to a neighborhood of the zero section 
of the associated bundle $G\times_{G_p}N_p$, where $N_p\subset T_p^*P$ is
the conormal space to $G.p$. Therefore, for some open neighborhood $V_p$
of $0\in N_p$, one has a homeomorphism $V_p/G_p\cong U_p$ onto some 
neighborhood of $x$ in $M$. The multiplicity of $M$ is simply the rank of 
a generic stabilizer for the $G$-action.   

Given any orbifold $M$, orbifold fiber bundles $\pi:E\ra M$ are defined 
by $H$-equivariant fiber bundles $Z\ra E_V\ra V$ 
in orbifold charts $(U,H,V)$, 
together with suitable compatibility conditions. 
Notice that the fibers $\pi^{-1}(x)$ are in general not diffeomorphic to 
$Z$, but only 
to some quotient of $Z$ by the action of the isotropy group 
$H_x$. For example, the tangent bundle of $M$ is an orbifold 
vector bundle, 
with fiber at $x\in M$ equal to $T_x V_x/H_x$, where $H_x$ is the 
isotropy 
group and $(U_x,V_x,H_x)$ an orbifold chart around $x$. 
Sections of an orbifold vector bundle $E$ are defined by $H$-invariant 
sections in orbifold charts. It is {\em not} true that any orbifold 
mapping $\sigma:\,M\ra E$ with $\pi\circ\sigma=\text{id}_M$ 
gives rise to a section 
of $E$. For example, $\C/\Z_2\ra\text{pt}$ does not have any 
non-vanishing sections.   

Continuing in this fashion, almost all constructions for manifolds 
can be generalized to orbifolds, 
by taking the corresponding $H$-equivariant version in local
orbifold charts. One thus defines suborbifolds, 
orbifold-principal bundles, 
de Rham theory, characteristic classes, Riemannian (complex, symplectic, 
Spin) structures, and so on. 
An orbifold principal bundle $H\ra P\ra M$ 
over $M$ is an orbifold $P$, together with a 
locally free action of some Lie group 
$H$, such that $M=P/H$. If $X$ is any $H$-space, one can form an 
associated 
bundle $P\times_H X=(P\times X)/H$. 
If only some finite cover $\hat{H}\ra H$ 
acts on $X$, one can still form the associated bundle, by merely 
regarding  
$P$ as a $\hat{H}$-orbifold principal bundle. Notice that if $H$ is 
connected, this does not depend on the choice of the cover. 

Suppose a compact, connected Lie group $G$ acts on $M$, in other words, 
$G\times M\ra M$ is an orbifold mapping. Since $G$ is connected, the 
components $F$ of the fixed point set $M^G$ are suborbifolds of $M$. 
If $G$ is abelian the normal bundle $\nu_F$ is even dimensional, and 
admits an invariant Hermitian structure (see e.g. \cite{BGV92}, p.217). 
If $M$ is oriented, 
any choice of such a Hermitian structure equips $F$ with an orientation.
 
Consider an orbifold chart $(U_x,H_x,V_x)$ around a point $x\in F$.
It is not always true that the $G$-action on $U_x$ lifts to $V_x$, 
but some finite covering $\hat{G}\ra G$ does. In particular, 
the weights for 
the action on $\nu_F(x)$ are weights for $\hat{G}$, but not 
necessarily for 
$G$. In \cite{LT94}, these are called ``orbi-weights'' of $G$.

Let $\A_G(M)$ be the space 
of equivariant differential forms, i.e. polynomial  
$G$-equivariant mappings $\alpha:\,\G\ra\A(M)$. Let $\d_\G$ be 
the equivariant 
differential
\begin{equation}\d_\G: \,\A_G(M)\ra \A_G(M), 
\,\d_\G(\alpha)(\xi)=\d\alpha(\xi)
-2\pi i\,\iota(\xi_M)\alpha(\xi),\end{equation}
where $\xi_M$ is the fundamental vector field corresponding to $\xi$.
Since $\d_\G^2=0$, $(\A_G(M),\d_\G)$ is a complex, and  
its cohomology $H_G(M)$ is called the {\em equivariant cohomology} 
of $M$. More generally, one can consider the complex $\A_G^\omega(M)$ 
of equivariant differential forms $\alpha$ which are defined and  
analytic around $0\in\G$, and define its cohomology $H_G^\omega(M)$.     

An example for an equivariantly closed form is the equivariant 
Euler form of the normal bundle $\nu_F$
. Suppose 
$G=T$ is abelian. If $\nu_F$ splits into a direct sum of invariant 
orbifold-line bundles $L_j$, with Chern classes $c_j$ and orbi-weights 
$\alpha_j$, one defines
\begin{equation}\Eul_\T(\nu_F,\xi)=\prod_j (c_j+2\pi i \l \alpha_j,\xi\r).\end{equation}
In the general case, this equation defines 
$\Eul_\T(\nu_F,\cdot)$ 
by means of the splitting principle.  

Suppose now that $M$ is compact, connected, and oriented, 
and consider the integration mapping 
$\int:\A_G^\omega(M)\ra \A_G^\omega(\text{pt})$. 
Since $\int\circ\,\d_\G=0$, this 
descends to a mapping on $H^\omega_G(M)$. Let $\iota_F:F\ra M$ 
be the embedding of the fixed point orbifolds.  

\begin{theorem}\label{orbiloc} 
(Orbifold version of the Berline-Vergne Localization Formula.)
Suppose $G=T$ is abelian, and let $\alpha\in \A_T^\omega(M)$ 
be $\d_\T$-closed. Then 
\begin{equation}\f{1}{d_M}\int_M \alpha(\xi)=\sum_F \f{1}{d_F}\int_F
\f{\iota_F^*\alpha(\xi)}
{\Eul_\T(\nu_F,\xi)}\label{locform}\end{equation}
for all $\xi$ in sufficiently small neighborhood of $0\in\T$. 
Here, the sum is over the fixed point orbifolds, 
and $d_F$ is the multiplicity of $F$ as a suborbifold of $M$. 
\end{theorem}
The proof of this result proceeds exactly as in the manifold case, see 
e.g. \cite{BGV92}. The mul\-ti\-pli\-ci\-ties 
occur because the proof involves an 
integration over the fibers of $\nu_F$, to be performed in local orbifold 
charts. In passing to the quotient, one has to divide by the ``twist''
of $\nu_F$, which gives a factor $d_M/d_F$.  

If $M$ is a $T$-manifold resp. orbifold with an (oriented) 
$T$-invariant        
boundary $Z=\partial M$, the localization formula  
includes additional boundary terms. Let us just consider the 
case $T=S^1$, and
suppose for simplicity that the action 
of $S^1$ on $\partial M$ is locally free. 
Let $j:\,Z\hra M$ be the inclusion, and $\theta\in\A^1(Z)$ 
a connection. 
Then the boundary contribution is given by 
\begin{equation}- \f{1}{d_M}\,\int_{\partial M}\,
\f{j^*\alpha(\xi)}{\d_\R\,\theta(\xi)}\,\wedge \theta.\end{equation}
For a proof (in the manifold case), see e.g. Kalkman \cite{K93}, 
or Vergne \cite{V94a}.  
\section{Riemann-Roch Theorems for Orbifolds}\label{RRO}

Suppose $M$ is a compact orbifold, equipped with a positive K\"ahler 
structure, and that $E\ra M$ is a holomorphic Hermitian orbifold vector
bundle over $M$. Just as in the manifold case, one has a twisted 
Dolbeault 
complex of $E$-valued forms,
\begin{equation}\bar{\partial}:\,\A^{i,j}(M,E)\ra \A^{i,j+1}(M,E),\end{equation}
and a corresponding $\Spin^c$-Dirac operator 
\begin{equation}\dirac_{\C}={\sqrt{2}}(\bar{\partial}+\bar{\partial}^*):
\,\A^{0,even}(M,E)\ra \A^{0.odd}(M,E).\end{equation}
If $M$ is only almost K\"ahler, it is still possible to 
construct $\dirac_{\C}$, 
using Hermitian connections on $E$ and on the canonical 
line bundle, $\wedge^{0,n}\,T^*M$. For the details of this construction 
see e.g. Duistermaat \cite{D95}. 
The dimension of the index space 
\begin{equation}\text{ind}(\dirac_\C)=[ker(\dirac_\C)]-[ker(\dirac_C^*)]\end{equation}
is independent of all choices, and will be called the Riemann-Roch number 
$\RR(M,E)$ of $E\ra M$.\\ 

\begin{remark}
If $M$ is a symplectic orbifold and $E\ra M$ a complex ( or
symplectic) vector bundle, one can always choose a compatible positive
almost complex structures on $M$ and Hermitian structure on $E$.  S
ince any two choices are homotopic, $\RR(M,E)$ does not depend on these 
choices. 
\end{remark}

If a compact Lie group $G$ acts on the above data and if all choices 
are made $G$-equivariant, one has $G$-actions on $\A^{i,j}(M,E)$, and 
the index space carries a virtual representation of $G$. In this case, 
we let  
$\RR(M,E)=\chi\in R(G)$ be the equivariant Riemann-Roch number, i.e. 
$\chi(g)=\text{tr}(g|\text{ind}(\dirac_\C))$. 
The orbifold version of the equivariant index theorem 
expresses $\chi$ as an integral 
of certain characteristic classes 
over an associated orbifold $\ti{M}$. 
Let us digress on how $\ti{M}$ is defined. As a set,  
\begin{equation}\ti{M}=\bigcup_{x\in M}\,\text{Conj}\,(H_x),\end{equation}
where $\text{Conj}(H_x)$ is the set of conjugacy classes in the 
isotropy group 
$H_x$. Given an atlas for $M$, orbifold charts for $\ti{M}$ can be 
constructed as follows. 
For each $V_i$, let 
\begin{equation}\hat{V}_i=\{(v,h)\in V_i\times H_i|\,h.v=v\}.\end{equation}
Then $H_i$ acts on $\hat{V}_i$ via $a.(v,h)=(a.v,\,a\,h\,a^{-1})$, and 
by definition $\ti{U}_i:=\hat{V}_i/H_i$. Notice that the 
action of $H_i$ on 
the preimage of a given connected component need not be effective, 
even if the action on $V_i$ was. 

One can show that the orbifold charts $(\ti{U}_i,\hat{V}_i,H_i)$ 
inherit the compatibility conditions
from the $(U_i,V_i,H_i)$, and thus define an orbifold 
structure on $\ti{M}$. 
Usually, $\ti{M}$ has various components of different dimension.

\begin{remark} (Properties of $\ti{M}$): \label{Holonomy}
\begin{enumerate}
\item
Let $G\times P\ra P$ be a locally free action of a compact Lie group 
$G$ on an or\-bi\-fold $P$. Then 
$M=P/G$ is an orbifold, and $\ti{M}=\hat{P}/G$, where 
\begin{equation}\hat{P}=\{(p,g)\in\ti{P}\times G|\,g.p=p\}.\label{hatP}
\end{equation}
\item
Consider the natural mapping $\tau:\ti{M}\ra M$. Since all mappings 
$\hat{V}\ra V$ in the above or\-bi\-fold charts are $H$-equivariant 
immersions (on all connected components), it follows that 
$\tau$ is an 
immersion. 
(If $M$ is a quotient of a manifold by a locally free action of an 
{\em abelian} group, $\tau$ restricts to orbifold embeddings 
of the connected components of $\ti{M}$. 
This is however false in general, because the 
orbifold isotropy groups are only defined up to isomorphism, and 
may be glued together globally in a nontrivial way.)

Let $N_{\ti{M}}\ra \ti{M}$ denote 
the normal bundle of this immersion.  
In a local orbifold chart $(U,V,H)$, $N_{\ti{M}}$ is obtained from
the normal bundle $N_{\hat{V}}$ of the immersion $\hat{V}\ra V$ 
by dividing out the action of $H$. 

\item 
Let $\ti{E}\ra\ti{M}$ be an orbifold vector bundle, given in local 
orbifold charts $(U,V,H)$ by $H$-equivariant vector bundles
$\ti{E}_{\hat{V}}\ra\hat{V}$. Observe that for all $(v,h)\in \hat{V}$, 
there is an action of $h$ on the fiber of $\ti{E}_{\hat{V}}$ 
over $(v,h)$. 
These actions glue together to give a canonical section $A$ of the 
automorphism bundle $\text{Aut}(\ti{E})$. (Notice however that the 
action of $A$ on sections of $\ti{E}$ is trivial!)
Suppose that $\ti{E}$ is a complex Hermitian orbifold bundle 
with connection.  
Let $F(\ti{E})\in\A^2(\ti{M},\,\text{End}(\ti{E}))$ be the curvature. 
We define 
{\em twisted} characteristic forms $\Ch^{\ti{M}}(\ti{E})$ 
and $D^{\ti{M}}(\ti{E})$ by 
\begin{equation}{\Ch}^{\ti{M}}(\ti{E})=\text
{tr}\big(A\,e^{\f{i}{2\pi}F(\ti{E})}\big)\in\A(\ti{M}).\end{equation}
and
\begin{equation}{D}^{\ti{M}}(\ti{E})=\det\big(1-A^{-1}\,e^{-\f{i}{2\pi}F(\ti{E})}\big)
\in\A(\ti{M}) .\end{equation}

\item
If $F$ is a connected suborbifold of $M$, the principal stratum of $F$ 
is contained in  
some stratum of $M$, and the multiplicity $d_F$ is the order of the  
isotropy group corresponding to that stratum. Moreover, the 
associated orbifold 
$\ti{F}$ is a suborbifold of $\ti{M}$. 
\end{enumerate}\end{remark}

Let us now turn to the case where $M$ is a symplectic orbifold. 
Then $\ti{M}$ is a symplectic orbifold, and 
the normal bundle $N_{\ti{M}}$ is a symplectic orbifold bundle. 
Choose a compatible positive almost complex structure $J$ on $M$, 
thereby making the 
tangent bundles of $M$ and $\ti{M}$ into Hermitian vector bundles. 
Consider a Hermitian connection on $TM$, with curvature 
$F(M)\in\A^2(M,\text{End}(TM))$,
and let 
\begin{equation}\Td({M})=\det\Big( \f{  \f{i}{2\pi}F(M)                        }
                {       1-e^{-\f{i}{2\pi}F(M)}          }\Big)
\end{equation}
and $\Td(\ti{M})$ be the corresponding Todd forms. 

Let $E\ra M$ be an orbifold vector bundle over an 
almost complex orbifold, and let $\ti{E}=\tau^* E$.
\begin{theorem}[Kawasaki]
The Riemann-Roch number $\RR(M,E):=\dim(\text{ind}(\dirac_\C))$ is 
given by the formula
\begin{equation}\RR(M,E)=\int_{\ti{M}}\f{1}{d_{\ti{M}}}\f{\Td(\ti{M})\,
{\Ch}^{\ti{M}}
(\ti{E})}{{D}^{\ti{M}}(N_{\ti{M}})}.\end{equation}
\end{theorem}
(For orbifolds $M$ that can be represented as quotients of manifolds
by locally free actions of {\em abelian} groups, the orbifold index 
theorem is due to Atiyah \cite{A74}.)
Let now $G$ be a compact, connected Lie group acting on $E\ra M$, 
and suppose all the above choices have been made $G$-invariant. By 
lifting the fundamental vector fields, one obtains an action on 
$\ti{M}$ of some finite cover of $G$ (due to the holonomy phenomenon 
mentioned in Remark \ref{Holonomy} (b), $G$ itself need not act on $\ti{M}$).
In 
\cite{V94}, M. Vergne has proved the following orbifold version of the 
equivariant index theorem of Atiyah-Segal-Singer \cite{AS68b}.
Denote by $\Td_\G(\ti{M},\xi)$, 
${\Ch}^{\ti{M}}_\G(\ti{E},\xi)$ etc. the 
equivariant characteristic classes for the respective bundles on 
$\ti{M}$. 
\begin{theorem}[Vergne] 
For $\xi\in\G$ sufficiently small, the equivariant Riemann-Roch number 
$\chi=\RR(M,E)$ is given by the formula 
\begin{equation}\chi(e^\xi)=\int_{\ti{M}}
\f{1}{d_{\ti{M}}}\f{\Td_\G(\ti{M},\xi)\,
{\Ch}^{\ti{M}}_\G
(\ti{E},\xi)}{{D}^{\ti{M}}_\G(N_{\ti{M}},\xi)}.\label{BV}\end{equation}
\end{theorem}
More generally, Vergne has proved a cohomological formula formula 
for $\chi(g\,e^\eta)$, with $\eta$ in the Lie algebra of the 
stabilizer of 
$g$, as an integral over $\tilde{M}^g$. If we take $\eta=0$, and 
$g=\exp(\xi)\in T$ generic in the sense that $g$ generates $T$, 
this includes the following special case:     
\begin{theorem}\label{ABF} (Fixed point formula for orbifolds.)
The equivariant index $\chi(e^\xi)$  of $\dirac_\C$ is 
given by the formula 
$\chi\,|T=\sum_F \,\chi_F$, 
where the sum is over the connected components $F$ of the 
fixed point orbifold $M^T$, and $\chi_F$ is a meromorphic function 
given by the formula 
\begin{equation}\chi_F(e^\xi)=\int_{\ti{F}}
\f{1}{d_{\ti{F}}} \f{ \Td(\ti{F})
\,{\Ch}^{\ti{F}}_\T(\ti{E}|\ti{F},\xi)}
{D^{\ti{F}}(N_{\ti{F}}) \, {D}_\T^{\ti{F}}(\ti{\nu}_F,\xi) },
\label{locc}
\end{equation}
Here $\nu_{F}$ is the normal bundle of $F$ in ${M}$, and   
$\ti{\nu}_F$ its pullback to $\ti{F}$. 
\end{theorem}
Notice that $\ti{\nu}_{F}$ is {\em not} the normal bundle of $\ti{F}$ 
in $\ti{M}$.
Observe also that Theorem \ref{ABF} follows from (\ref{BV}), by 
applying the localization theorem.\\

\noindent{\bf Example:}\\
The teardrop orbifold $M$ of order $k$ is obtained by gluing a copy 
of $\C/\Z_k$, with coordinate $w^k$, with a copy of $\C$, with 
coordinate $z$, via the mapping $z\mapsto w^{-k}$ for $z\not=0$. 
Note that $\ti{M}$ is a disjoint union of a copy of $M$ and 
$(k-1)$ points. 
Let $G=S^1$ act by rotation, i.e. $z\mapsto e^{i\phi}\cdot z$, and let 
$E=M\times \C$ be the trivial line bundle. The fixed point set 
consist of two points, $z=0,\infty$.
At $z=0$, $M$ is smooth, the character of the action on the normal bundle 
is $+1$, and so the fixed point contribution becomes 
$\chi_{\{0\}}(e^{i\phi})=(1- e^{i\phi})^{-1}.$
At $z=\infty$, there is an orbifold singularity of order $k$, hence
$\ti{F}$ consists of $k$ points. The normal bundle 
is isomorphic to $\C/\Z_k$, and $S^1$ acts by the orbi-weight 
$-\f{1}{k}$. 
Thus 
$$  \chi_{\{\infty\}}(e^{i\phi})=\f{1}{k}\sum_{l=0}^{k-1} 
\f{1}{1-c^{l}\,e^{-\f{i}{k}\phi}},$$ 
where $c=\exp(-\f{2\pi i}{k})$. Using the identity   
\begin{equation}\f{1}{k}\sum_{l=0}^{k-1} \f{1}{1-c^{l}\,u}=\f{1}{1-u^k},
\label{magic}\end{equation}
we can carry out the summation and find 
$\chi_{\{\infty\}}(e^{i\phi})=(1-e^{-i\phi})^{-1}$.  The 
two contributions add up to give $\chi(e^{i\phi})=
(1-e^{i\phi})^{-1}+(1-e^{-i\phi})^{-1}=1$. 
\\[.2in] We will now give a simple application of the fixed point
formula in connection with holomorphic induction.  Let $G$ be a
compact connected Lie group, with maximal torus $T$, and
$\T^*_+\subset \T^*=(\G^*)^T$ some choice of a positive Weyl chamber.
For each face $\Sigma$ of $\T^*_+$, there is a compact, connected
subgroup $G_\Sigma\subset G$ with the property that
$G_\alpha=G_\Sigma$ for all $\alpha\in \text{int}(\Sigma)$.  Write
$\T^*_{\Sigma,+}\supset \T^*_+$ for the positive Weyl chamber of
$G_\Sigma$, $W_\Sigma\subset W$ for the Weyl group,
$\lie{R}_{\Sigma,+}\subset \lie{R}_{+}$ for the positive roots, and
$\Lambda^*_{\Sigma,+}=\Lambda^*\cap \T^*_{\Sigma,+}$ for the dominant
weights.  
Let $G/G_\Sig$ be equipped with its canonical
$G_\Sig$-invariant complex structure, corresponding to the
interpretation as a coadjoint orbit for $G$. If $Y_\Sigma$ is a
compact, almost complex $G_\Sig$-orbifold, the associated bundle
$G\times_{G_\Sig}Y_\Sigma$ has a canonically induced almost complex
structure, and given a $G_\Sig$-equivariant orbifold-vector bundle
$E_\Sigma\ra Y_\Sigma$ one can form the associated bundle
$G\times_{G_\Sig} E_\Sigma\ra G\times_{G_\Sig}Y_\Sigma$.
In particular, one can apply this to the case where
$Y_\Sig=G_\Sigma\cdot\mu$ is an integral coadjoint orbit through a
$G_\Sig$-weight $\mu\in \Lambda^*_{\Sigma,+}$, and
$E_\Sigma=G_\Sigma\times_T\,\C_\mu$ the corresponding pre-quantum line
bundle.  Then $\RR(G_\Sigma\cdot\mu,\,G_\Sigma\times_{G_\mu}\,\C_\mu)$
is the irreducible $G_\Sigma$-representation $\chi_{\Sigma,\mu}$ 
labelled by $\mu$. The associated
bundle $G\times_{G_\Sigma}Y_\Sigma$ is just the coadjoint orbit
$G\cdot\mu$, and the Riemann-Roch number of 
$G\times_{G_\Sigma}E_\Sigma=G\times_{G_\mu}\C_\mu$ is by definition the 
holomorphic induction of $\chi_{\Sigma,\mu}$. 
Let 
$$ \Ind_{G_\Sigma}^G:\,R(G_\Sigma)\ra R(G) $$
denote the holomorphic induction map.

\begin{theorem}\label{slice}
Let $Y_\Sigma$ be a compact almost complex $G_\Sigma$-orbifold, amd 
$E_\Sigma\to Y_\Sigma$ a $G_\Sigma$-equivariant orbifold vector bundle. 
The $G$-equivariant Riemann-Roch number of $G\times_{G_\Sigma}E_\Sigma$ is related to the $G_\Sigma$-equivariant Riemann-Roch 
number of $E_\Sigma$ by holomorphic induction:
\begin{equation}
\RR(G\times_{G_\Sigma}Y_\Sigma,\, 
G\times_{G_\Sigma}\,E_\Sigma)=\Ind_{G_\Sigma}^G\,\,\RR(Y_\Sigma,E_\Sigma)
\label{slicef}
\end{equation}
\end{theorem}
\begin{proof}
The proof is by applying the fixed point formula to both sides.
We write $\chi= \RR(G\times_{G_\Sigma}Y_\Sigma,\, 
G\times_{G_\Sigma}\,E_\Sigma)$ and $\chi_\Sig=\RR(Y_\Sigma,E_\Sigma)$, 
and denote by $N_\Sig(\mu)$ the multiplicity of a $G_\Sigma$-weight 
$\mu\in\Lambda^*_{\Sigma,+}$ in $\chi_\Sigma$.  
Since $(G/G_\Sig)^T=W/W_\Sig$, 
$$ M^T=(G\times_{G_\Sig} Y_\Sigma)^T=W\times_{W_\Sigma}Y_\Sigma^T.$$
We will think of $W/W_\Sig$ as the set of all $w\in W$ such that 
$w\cdot\T^*_+\subset \T^*_{\Sig,+}$, so that 
$$ \chi(e^\xi)=\sum_{F\subset M^T}\chi_F(\xi)=\sum_{w\in W/W_\Sig}
\sum_{F\subset Y_\Sigma^T}\chi_{F}(w^{-1}(\xi)).$$
Consider the $T$-action on $T_xM$, for $x\in Y_\Sigma^T\subset
M^T$. The weights for the action on $T_xM/T_xY_\Sigma$ are just
$\lie{R}_+ - \lie{R}_{\Sigma,+}$, the set of positive roots of $G$
that are not positive roots for $G_\Sigma$.  It follows that for all
$F\subset Y_\Sigma^T\subset M^T$, the fixed point contributions
$\chi_{\Sig,F}(\xi),\,\chi_F(\xi)$ with respect to $Y_\Sigma,\,M$ are
related by
$$ \chi_F(\xi)=\f{\chi_{\Sig,F}(\xi)}{
\prod_{\beta\in\lie{R}_+ 
-\lie{R}_{\Sigma,+}} (1-e^{-2\pi i\l\beta,\xi\r})}.$$
Summation over all fixed point contributions in $Y_\Sigma^T\subset M^T$ gives, 
by another application of the fixed point formula, 
\begin{eqnarray*}
\chi(e^\xi)&=&\sum_{w\in W/W_\Sig} \sum_{F\subset
Y_\Sigma^T}\f{\chi_{\Sig,F}(\xi)}{ \prod_{\beta\in\lie{R}_+
-\lie{R}_{\Sigma,+}} (1-e^{-2\pi i\l\beta,\xi\r})}\\ &=&\sum_{w\in
W/W_\Sig}\f{\chi_\Sigma(e^{w^{-1}(\xi)})}{\prod_{\beta\in\lie{R}_+
-\lie{R}_{\Sigma,+}} (1-e^{-2\pi i\l\beta,\xi\r})}\\
&=&
\sum_{\mu\in\Lambda^*_{\Sig,+}}\,N_{\Sigma}(\mu)\sum_{w\in
W/W_\Sig} \f{\rho_{\Sig,\mu}
(e^{w^{-1}(\xi)})}
{\prod_{\beta\in\lie{R}_+ 
-\lie{R}_{\Sigma,+}} (1-e^{-2\pi i\l\beta,w^{-1}\xi\r})}.
\end{eqnarray*}
The same formula with $N_{\Sig}(\mu)=1$ gives an expression for 
$\Ind_{G_\Sig}^G\,\rho_{\Sig,\mu}$. Thus 
$$\chi(e^\xi)=\sum_{\mu\in\Lambda^*_{\Sig,+}}\,N_{\Sigma}(\mu) 
\Ind_{G_\Sig}^G\,\rho_{\Sig,\mu}(e^\xi)
=\Ind_{G_\Sig}^G\, \chi_{\Sig}(e^\xi).
$$
\end{proof}

\section{Symplectic Surgery}\label{SymplecticSurgery}

Let $(M,\omega)$ be a symplectic orbifold, and $S^1\times M\ra M$ a 
Hamiltonian action with moment map $\Phi:\,M\ra \R$. 
Suppose that 0 
is a regular value of $\Phi$, and let $M_{red}=\Phi^{-1}(0)/S^1$
be the symplectic quotient. 
The symplectic cutting construction 
of Lerman \cite{L93} yields two new symplectic orbifolds $M_-,\,M_+$, 
such that, as a set, 
\begin{equation}M_-=M_{red}\cup\{\Phi<0\},\,\,M_+=
M_{red}\cup\{\Phi>0\},\end{equation}
and 
the embeddings of $M_{red},\,\{\Phi<0\},\,\{\Phi>0\}$ are symplectic.  

The structure of symplectic orbifolds on these spaces is obtained 
as follows. 
Consider the product $M\times\C$, with symplectic form 
$\omega-\f{i}{2}\d z\wedge\d\bar{z}$, and the diagonal circle action 
$e^{i\phi}\cdot(x,z)=(e^{i\phi}\cdot\,x,\,e^{-i\phi}z)$,
with moment map $\psi(x,z)=\Phi(x)-|z|^2$. Then $0$ is a regular value of 
$\psi$, and so $\psi^{-1}(0)/S^1$ is a symplectic orbifold. 
The level set $\psi^{-1}(0)$ consists of two components:
\begin{equation}\label{TwoComp}
\psi^{-1}(0)=\{(x,z)|\,\Phi(x)>0,\,\,|z|^2=\Phi(x)\}\cup\,
(\Phi^{-1}(0)\times\{0\}).\end{equation}
It is now easy to see that the map
\begin{equation}\alpha:\{\Phi>0\} \to \psi^{-1}(0)\subset M\times \C,\,x
\mapsto (x,\,\sqrt{\Phi(x)})\label{alpha}\end{equation}
satisfies $\alpha^*(\omega-\f{i}{2}\d z\wedge\d\bar{z})=\omega$.  
This identifies $\psi^{-1}(0)/S^1=M_+$. For $M_-$, one simply takes the 
opposite circle action on $\C$.\\ 

\noindent{\bf Example}: Take $M=\C P(1)$, with the Fubini-Study form, 
and consider the natural action of $S^1$, $e^{i\phi}.[z_0:\,z_1]=
[e^{i\phi}z_0:\,z_1]$. A moment map for this action is given by 
$\Phi([z_0:z_1])=-\f{|z_0|^2}{||z||^2}+\f{1}{2}$. Zero is a regular 
value of $\Phi$, the action of $S^1$ on $\Phi^{-1}(0)$ is free, and 
the cut spaces are spheres with half the volume of $\C P(1)$. 
Consider, on the other hand, the moment map $\Phi^{(k)}=k\Phi$, 
for $k\in \N$, which lets $S^1$ rotate $\C P(1)$ with $k$-fold speed. 
For this case, the $S^1$ action on 
the zero level set is only locally free, and the cut spaces are 
teardrop-orbifolds with a $\Z/\,k\Z$-singularity at $M_{red}=\{pt.\}$.\\

For the cutting construction, 
the Hamiltonian $S^1$ action needs only be locally defined, in some 
neighborhood of the hypersurface $Z=\Phi^{-1}(0)$. The corresponding 
cut space $M_{cut}$ (equal to $M_+\cup M_-$ for a globally defined  
$S^1$ action) 
has one or two connected components, depending on 
whether or not $Z$ disconnects
$M$. We denote the two copies of $M_{red}$ in $M_{cut}$ by 
$M_{red}^\pm$. 

Consider for example the torus  $M=T^2=\R^2/\Z^2$, with the action of 
$S^1$ defined by $e^{2\pi i t}.[x,y]=[x+t,y]$. In a neighborhood
of the circle $e^{2\pi i t}.[0,0]$, one can take $\Phi([x,y])=y$ as 
a moment 
map, and $M_{cut}$ is symplectomorphic to $S^2$.\\

\begin{remark}
\begin{enumerate}
\item
If a Lie-group $G$ acts on $M$ in a symplectic fashion, 
and if this action commutes with the action of $S^1$, one obtains a 
symplectic 
$G$-action on $M_{cut}$. If the action on $M$ has a moment map, the 
restriction of this moment map to $M-Z$ extends 
smoothly to a $G$-moment map on $M_{cut}$. 
\item 
Consider $Z\ra M_{red}$ as an orbifold principal $S^1$-bundle.
The normal orbifold-bundle of $M_{red}^\pm$ in $M_{cut}$ 
is equal to the associated bundle $Z\times_{S^1}\C$, where $S^1$ 
acts on $\C$ by the character $e^{\mp i\phi}$. In particular, 
they have opposite 
Chern classes.
\item 
The cutting construction does not use nondegeneracy of $\omega$.
\end{enumerate}
\end{remark}

Let us return to the case where the $S^1$ action is globally defined, and 
consider an $S^1$-equivariant Hermitian orbifold
vector bundle
 $E\ra M$. Let $\pi_1:\,M\times\C \ra M$ denote projection to the 
first factor, and define $E_{red}=(E|\Phi^{-1}(0))/S^1$ and 
$E_{\pm}=(\pi_1^*\,E|\,\psi^{-1}(0))/S^1$, with the induced Hermitian 
structure. The mapping (\ref{alpha}) induces a mapping 
$$  \,E|\,\{\Phi>0\}\times S^1\,\ra \pi_1^*\,E,$$ 
which descends to an isomorphism of Hermitian vector 
bundles $E_+|\{\Phi >0\}\cong E|\{\Phi >0\}$, and similarly of course  
$E_-|\{\Phi<0\}\cong E|\{\Phi <0\}$. It follows that cutting is local 
on the level of Hermitian vector bundles, in the following sense: 
Suppose that the $S^1$ action on $E\ra M$ is only defined near a given 
$S^1$-invariant hypersurface $Z=\Phi^{-1}(0)$. 
Then one can define 
a reduced Hermitian bundle $E_{red}\ra M_{red}$ and a cut bundle, 
$E_{cut}\ra M_{cut}$. If a compact Lie group $G$ acts on 
$E\ra M$ and if all 
data are $G$-invariant, one obtains smooth $G$-actions on 
$E_{red}$ and $E_{cut}$.   
Given this situation, let us now study the relation of the 
($G$-equivariant) 
Riemann-Roch numbers of $M$ and its cut- and reduced spaces.

\begin{theorem}
[Gluing Formula] 
The $G$-equivariant Riemann-Roch numbers of $E$, $E_{cut}$, 
$E_{red}$ are related by   
 \begin{equation}\RR(M,E)=\RR(M_{cut},E_{cut})-\RR(M_{red},E_{red}).\label{glue1}\end{equation}
\end{theorem}
\begin{proof}
Let us assume $G=\{1\}$ for simplicity of notation; for the general 
case, one simply replaces characteristic classes by $G$-equivariant
characteristic classes. Let $B\subset M$ be a tubular 
neighborhood of $Z$, 
such that the Hamiltonian $S^1$ action 
is defined and locally free on $B$.
Pick compatible almost complex structures of $M$ and $M_{cut}$ 
in such a way 
that they agree over some neighborhood of $M-B$. Also, choose Hermitian 
connections $\nabla,\,\nabla_{cut}$ on $E,\,E_{cut}$ that agree 
over some 
neighborhood of $M-B$ and are $S^1$-invariant over $B$. 

Let $\ti{B}$ be the preimage of $B$ in $\ti{M}$. 
Over $\ti{B}$, we can introduce $S^1$-equivariant 
characteristic classes, and write
$$ \RR(M,E)=\int_{\ti{M}-\ti{B}}\f{1}{d_{\ti{M}}}
\f{\Td(\ti{M})\,\Ch^{\ti{M}}(\ti{E})}{D^{\ti{M}}(N_{\ti{M}})}
   +
\int_{\ti{B}}\f{1}{d_{\ti{M}}}
\f{\Td_{\R}(\ti{M},\xi)\,\Ch^{\ti{M}}_{\R}(\ti{E},\xi)}
{D_{\R}^{\ti{M}}(N_{\ti{M}},\xi)}\Big|_{\xi=0}.$$ 
To the second term, we apply the localization formula for manifolds 
with boundary, and rewrite it in the form 
$$  \int_{\p B}\alpha(\xi)\,\Big|_{\xi=0},$$ 
for some equivariant differential form $\alpha\in\A_\R(\p B)$. 
Let us do the same computation for $M_{cut}$, with $B$ 
replaced by $B_{cut}$. 
$M_{cut}$ contains two copies of $M_{red}$ as codimension 2, 
$S^1$-fixed suborbifolds, 
but with opposite normal bundles, and the complement is symplectomorphic 
to $M-Z$. We hence obtain the same terms as above, plus two additional 
terms corresponding to the two copies of $M_{red}$:
\begin{equation}\int_{\ti{M}_{red}}\f{1}{d_{\ti{M}_{red}}}
\f{\Td(\ti{M}_{red})\,
\Ch^{\ti{M}_{red}}_{\R}(\ti{E}_{red},\xi)}
{D_{\R}^{\ti{M}_{red}}(N_{\ti{M}_{red}})\,D_\R^{\ti{M}_{red}}
(\ti{\nu},\xi)}
+
\int_{\ti{M}_{red}}\f{1}{d_{\ti{M}_{red}}}
\f{\Td(\ti{M}_{red})\,\Ch^{\ti{M}_{red}}_{\R}
(\ti{E}_{red},\xi)}
{D_{\R}^{\ti{M}_{red}}(N_{\ti{M}_{red}})\,
D_\R^{\ti{M}_{red}}(\ti{\nu}^*,\xi)}
,\end{equation}
evaluated at $\xi=0$. But  
\begin{equation} D_\R^{\ti{M}_{red}}(\ti{\nu},\xi)^{-1}+ 
D_\R^{\ti{M}_{red}}(\ti{\nu}^*,\xi)^{-1}=1,
\label{loccont}\end{equation}
(using the formal identity $(1-z)^{-1}+(1-z^{-1})^{-1}=1$), so
the sum of these terms is just $\chi_{red}$. 
\end{proof}
 
\begin{remark} 
\begin{enumerate}
\item
Notice that in the above proof, the Chern character 
may be replaced by any characteristic class of $\ti{E}$. 
\item
For any symplectic  orbifold $M$, let 
$\tau(M)=\int_M \Td(M)$. 
The above Theorem, applied to the trivial bundle $E=M\times\C$, 
sometimes greatly simplifies the computation of $\tau(M)$. 
\begin{enumerate}
\item
Consider for example
the cutting of  $M=\C P(n)$ along the equator sphere $Z=S^{2n-1}$. 
Then $M_{cut}$ 
consists of two copies of $\C P(n)$, and the reduced space is 
$\C P(n-1)$. 
Hence $\tau(\C P(n))=2 \tau(\C P(n))-\tau(\C P(n-1))$, or $\tau(\C P(n))=
\tau(\C P(n-1))=\ldots = \tau(\C P(0))=1$. 
\item 
Consider next a Riemann surface $M$ of genus $g$. 
By cutting along $g$ circles, $M$ can be made into a sphere, 
hence $\tau(M)=\tau(S^2)-g=1-g$.  
\item
Let $M$ be any symplectic manifold, and $p\in M$. 
The symplectic analogue of the blow-up $Bl_p(M)$ of $M$ at $p$
(see \cite{GS89}) 
replaces $p$ by a ``small'' $\C P(n-1)$, and 
may be regarded as a cutting operation (see \cite{L93}), 
with $M_+=Bl_p(M)$, $M_-=\C P(n)$, and $M_{red}=\C P(n-1)$.  
Since $\tau(\C P(k))$ is equal to one, 
$\tau(Bl_p(M))=\tau(M)-\tau(\C P(n))+\tau(\C P(n-1))=
\tau (M)$. 
More generally, one can consider symplectic blow-ups $Bl_N(M)$ along   
a symplectic submanifold $N\subset M$, and a similar argument 
shows that $\tau(Bl_N(M))=\tau(M)$.
\item
Let $M_k$ be the teardrop-orbifold of order $k$. Recall from the 
first example in this section that $M_k$ can be constructed by 
cutting the sphere $\C P(1)=S^2$, with $S^1$ acting with $k$-fold 
speed. The gluing formula $\RR(\C P(1),\C)=2\,\RR(M_k,\C)-
\RR({\text pt},\C)$ shows (using Kawasaki's formula) that
\begin{eqnarray*} 
1&=&2\Big(\tau(M_k)+\f{1}{k}\sum_{j=1}^{k-1}
\f{1}{1-e^{2\pi i\f{j}{k}}}\Big)-1\\
&=&2\big(\tau(M_k)+\f{1}{2}(1-\f{1}{k})\big)-1.
\end{eqnarray*}
Hence $\tau(M_k)=\f{1}{2}\big(1+\f{1}{k}\big)$. More generally, if 
$M$ is a 2-dimensional symplectic orbifold of genus $g$, 
with $r$ orbifold-singularities of order $k_1,\ldots,k_r$, one finds 
\begin{equation}\tau(M)=1-g-\f{1}{2}\sum_{i=1}^r(1-\f{1}{k_i}).\end{equation}
\end{enumerate}
\item
As pointed out in \cite{L93}, the original space can be recovered 
from the cut space by means of the symplectic gluing procedure of 
Gompf \cite{G93}:
Let $N$ be a symplectic orbifold, and let $\iota_i:\,N\ra M_i$ be 
symplectic embeddings into two given symplectic orbifolds with 
$\dim M_i=\dim N+2$. If the corresponding normal bundles $\nu_i$ 
are opposite, $\nu_1^*\cong\nu_2$,
Gompf's method yields a new symplectic manifold $M$, 
which as a set is a union of $M_1-N,\,M_2-N$ and the unit circle 
bundle of $\nu_1=\nu_2^*$.
One can show that the gluing procedure is 
quantizable as well, i.e. 
given symplectic orbifold vector bundles $E_i$ with the 
same restriction to 
$N\subset M_i$, then there is a symplectic orbifold 
vector bundle $E$ on $M$, 
and the $E_i$ are obtained from $E$ by cutting.
\end{enumerate}   
\end{remark}
Let $\Z_l$ be the generic stabilizer for the action of $S^1$ 
on $\Phi^{-1}(0)$. Then the action of $H=S^1/\Z_l$ on $\Phi^{-1}(0)$ 
is generically free. We may use $H$ instead of the original $S^1$
for the cutting construction, and obtain a cut space $M_+$ which is 
again the disjoint union of $M_{red}$ and $\{\Phi>0\}$, but they 
are glued together in a different way. Conversely, starting from 
a generically free action on $\Phi^{-1}(0)$, with corresponding 
cut space $M_+$, obtains new cut spaces $M_+^{(l)}$ ($l\in\N)$ 
by replacing $S^1$ by its $l$-fold cover. The normal 
bundle $\nu^{(l)}$ of $M_{red}$ 
in $M_+^{(l)}$ is equal to the quotient of $\nu$ by the action 
of $\Z_l$, and the equivariant Chern class (where we 
consider the action of 
the {\em original} $G=S^1$ on $M_+^{(l)}$) gets divided by a factor 
of $l$. The space $\ti{M}_{red}^{(l)}$ consists simply of $l$ 
copies of $\ti{M}_{red}$, in particular $d_{\ti{M}_{red}^{(l)}}=
l\,d_{\ti{M}_{red}}$. Therefore, the fixed point contribution 
$\chi_{{M}_{red}}^{(l)}$ becomes 
$$  \chi_{{M}_{red}}^{(l)}(e^\xi)=  
\f{1}{l} \int_{\ti{M}_{red}}\f{1}{d_{\ti{M}_{red}}}
\f{\Td_{\R}(\ti{M_{red}})\,
\Ch^{\ti{M}_{red}}_{\R}(\ti{E}_{red},\xi)}
{D_{\R}^{\ti{M}_{red}}(N_{\ti{M}_{red}})\,}
\,\sum_{j=0}^{l-1} \f{1}{D_\R^{\ti{M}_{red,\,j}}(\ti{\nu}^{(l)},\xi)},$$ 
where $\ti{M}_{red,\,j}$ denotes the $j$th copy of $\ti{M}_{red}$.
The sum over $j$ can be computed, using Equation (\ref{magic}):
$$  \f{1}{l}\sum_{j=0}^{l-1} \f{1}{D_\R^{\ti{M}_{red,j}}
(\ti{\nu}^{(l)},\xi)}
=\f{1}{D_\R^{\ti{M}_{red}}(\ti{\nu},\xi)}.$$ 
This leads to the following  observation: 
\begin{proposition}
\label{observation}
The fixed point contributions $\chi_{M_{red}}(e^\xi)$ 
to $\RR(M_\pm,E_\pm)$ 
depend only  the cutting hypersurface 
$Z$, not on how the projection $Z\ra M_{red}$ is made into an 
orbifold $S^1$
bundle.
\end{proposition}

We will now consider the special case that $E=L$ is a 
pre-quantum line bundle, 
with Hermitian connection satisfying the pre-quantum condition 
$$ \f{i}{2\pi}\text{curv}(\nabla)=\omega,$$ 
and such that the $S^1$-action on $M$ lifts to $L$ according to 
(\ref{lift}). It is well known that $L_{red}=(L|\Phi^{-1}(0))/S^1$ 
has a unique Hermitian connection, such that the pullback to 
$\Phi^{-1}(0)$ is equal to the restriction of $\nabla$.  

The  trivial line bundle $L_\C =\C \times \C$, with Hermitian 
fiber metric 
$$  \l w_1,w_2\r=\exp(-\pi\,|z|^2)\,\bar{w}_1\cdot w_2,\,\,$$ 
($w_i\in (L_\C)_z$) and connection 1-form $A=-\pi\,z\d \bar{z}$ is 
a pre-quantum bundle for the symplectic structure 
$-\f{i}{2}\d z\wedge \d\bar{z}$. The pre-quantum lift (\ref{lift}) of the 
$S^1$ action on $\C $ to $L_\C$ is simply the trivial action. 
Hence, $L\boxtimes L_\C$ is a pre-quantum line bundle for 
$M\times \C $, 
and the above prescription gives a pre-quantum bundle for $M_+$, 
$$  L_+\,=\,(L\boxtimes L_\C|\psi^{-1}(0))/S^1.$$ 
Notice that over $\{\Phi>0\}$, we have now two pre-quantum line bundles: 
One coming from the embedding into $M$, the other from the embedding 
into $M_+$. 
\begin{theorem}
There exists a canonical isomorphism of pre-quantum line bundles 
with connection, 
$L_+|\{\Phi>0\}\cong L|\{\Phi>0\}$, and $L_+|M^+_{red}\cong L_{red}$.
\end{theorem}
\begin{proof}
Since 
$L_\C $ is trivial, $L\boxtimes L_\C \cong \pi_1^*\,L$ as a 
complex vector bundle,
but the fiber metric is multiplied by a factor $e^{-\pi\,|z|^2}$, and 
the connection is 
$$ \nabla\boxtimes \nabla_\C=\pi_1^*\nabla-\,\pi\,z\d \bar{z}.$$ 
The map  
\begin{equation}\alpha:\{\Phi>0\} \to \psi^{-1}(0)\subset M\times \C,\,x
\mapsto (x,\,\sqrt{\Phi(x)})\end{equation}
in (\ref{alpha}) is covered by a bundle map 
$$  
\beta:\,L|\{\Phi>0\}\to \pi_1^*L,\,\,
\lambda\mapsto \lambda \,e^{\f{\pi}{2}\Phi(x)}$$ 
which preserves the fiber metric, and satisfies 
$$  \beta^*(\nabla)=\nabla+\f{\pi}{2}\d\Phi,\,\,\,\,\beta^*(\pi z \d\bar{z})=\f{\pi}{2}\d\Phi.$$ 
This proves the first isomorphism, and the second is obvious.
\end{proof}
\noindent 
A similar result holds of course for $L_-\ra M_-$.  As a consequence
of this Theorem, cutting is local even on the level of pre-quantum
line bundles.

It is possible to cut at nonzero levels $\alpha\in\R$ of the 
moment map, by simply redefining the moment map used for the cutting
construction as $\Phi'=\Phi-\alpha$. 
Given a $S^1$-invariant complex vector bundle $E\ra M$, one again 
obtains  
cut bundles $E_+$ and $E_-$. Notice however that if $L$ is a pre-quantum 
line bundle for $M$, the cut bundles obtained in this way
are not pre-quantum line bundles for $M_\pm$. Similarly,  
$L|\Phi^{-1}(\alpha)/S^1$ is not a pre-quantum bundle for 
$M_\alpha=\Phi^{-1}(\alpha)/S^1$, 
since the $S^1$ action on $L$ satisfies the 
pre-quantum condition 
with respect to $\Phi$, but not with respect to $\Phi'$. 
If $\alpha$ is an integer, the correct lift is obtained by 
tensoring $L$ with the trivial line bundle $\C$, with $S^1$ acting 
by the character $e^{-2\pi \,i\l\alpha,\xi\r}$, before applying the 
cutting construction. More generally, if 
$\alpha$ is a rational number, one can choose $k\in N$ such 
that $k\alpha$ is an integer, and cut with respect to the $k$-fold cover
of the original $S^1$. This will introduce extra orbifold singularities 
in $M_\pm$, thereby making $M_\pm$ quantizable by means of orbifold-line 
bundles. 
 
To conclude this section, let us briefly explain (following 
\cite{DGMW95}) why the gluing 
formula implies the Guillemin-Sternberg conjecture for the case 
$G=S^1$ (and therefore also for the abelian case $G=T$, using 
reduction in stages). 
It can be deduced from the fixed point 
formula that the multiplicity of $0$ is equal to $\RR(M_0,L_0)$ if 
zero 
is a maximum or minimum of $\Phi$, since in that case $M_0$ is 
contained in $M$ as a fixed point manifold, and $L_0=L|M_0$.
Let $N_\pm$ be the multiplicity functions for the cut bundles 
$L_\pm\ra M_\pm$. 
By the gluing formula, 
$$  N(\mu)=N_+(\mu)+N_-(\mu)-\RR(M_0,L_0)\,\delta_{\mu,0}.$$ 
Since $0$ is a maximum of the moment map for $M_-$ and a 
minimum for $M_+$, we have $N_+(0)=N_-(0)=\RR(M_0,L_0)$. Hence 
$N(0)=\RR(M_0,L_0)$, q.e.d. In particular, a weight does not occur 
if it is not in the image of the moment map. 

\begin{corollary}
(See \cite{DGMW95})\label{corcut} 
Let $(M,L)$ be a quantizable Hamiltonian $S^1$ space, with 
multiplicity function $N(\mu)$, and suppose $0$ is a regular 
value of the moment map. Then the multiplicity function for 
the cut bundle $L_+\ra M_+$ is given by
$N_+(\mu)=N(\mu)$ if $\mu\ge 0$, and $N_+(\mu)=0$ if $\mu<0$.
\end{corollary}
 
\section{Multiple Cutting}

\subsection{Cutting with respect to polytopes}
 
Let $T$ be a $k$-torus, and 
$M$ a compact connected Hamiltonian $T$-orbifold, with moment 
map $\Phi:\,M\ra \T^*$. Denote 
\begin{equation}M_{\H}=\{x\in M|\,\T_x=\H\}.\end{equation}
The connected components of $M_\H$ are symplectic suborbifolds  
of $M$
\footnote{Given an action of a compact Lie group $G$ on an 
orbifold $M$, the set $M_{H}$ of points with stabilizer $G_x$ equal to 
some fixed group $H$ is in general not a suborbifold, because 
the action of $G_x$ 
does not always lift to an action in orbifold charts around $x$.}.  
Let $M=\cup M_i$ be the corresponding decomposition of $M$, where 
the $M_i$ are the connected components of the various $M_{\H}$, and let 
$\H_i$ denote the isotropy group corresponding to $M_i$. 

By the convexity theorem of Atiyah \cite{A82} and Guillemin-Sternberg 
\cite{GS82c}, and its extension 
to the orbifold case by Lerman and Tolman \cite{LT94}, the image 
$\Delta=\Phi(M)$ is equal to the convex hull of the image of the fixed 
point set, $\Phi(M^T)$. More precisely, the closure of $\Phi(M_i)$ 
is a convex polytope $\Delta_i$ of codimension $\dim \H_i$, 
which is contained in an affine 
space of the form $\{\alpha\}+\H_i^0\subset \T^*$. It turns out that 
the $\Phi(M_i)$ define a subdivision of $\Delta$ into rational convex 
polytopes. In particular, the  connected  components of the set 
$\Delta_*\subset \Delta $ of regular values are convex open 
subpolytopes of $\Delta$.

Let $R$ be a rational, convex polytope $R\subset \T^*$, 
i.e. a polytope defined by a finite number of inequalities and 
equalities, 
\begin{equation}\left\{\begin{array}{l@{\quad:\quad}l}
\l \alpha,v_i\r\ge \mu_i&i=1,\ldots,r\\
\l \alpha,v_i\r= \mu_i&  i=r+1,\ldots,N
\end{array}\right.\,,
\label{de}\end{equation}
where $v_i\in \Lambda$, and $\mu_i\in \R$. 
Denote by $R_\R$ the affine subspace generated by $R$, and 
by $T_{R}\subset T$ the torus perpendicular to $R_\R$.  
We will call $R$  {\em admissible}, if it satisfies the following 
conditions: 
\begin{itemize}
\item[(A)] 
The affine hyperplanes $\l \alpha,v_i\r= \mu_i$, $i=1,\ldots,N$, 
are all transversal.
\item[(B)] 
The faces of $R$ are transversal to all (interior and 
exterior) faces of $\Delta$.
\end{itemize}  
Let $\alpha\in R\cap\Delta$, and $x\in\Phi^{-1}(\alpha)$. Let 
$R_\alpha\subset R$ be the unique face that contains $\alpha$ 
in its interior. Condition (B) is equivalent to 
$\T_x^0+\T_{R_\alpha}^0=\T^*$, 
or to 
\begin{itemize}
\item[(B')] For all $\alpha\in R\cap\Delta$ and $x\in\Phi^{-1}(\alpha)$,
$ \T_x\cap \T_{R_\alpha}=\{0\}$.
\end{itemize}  
For the pre-quantization setting, we will need the additional 
condition 
\begin{itemize}
\item[(Q)] For all $i=1,\ldots,N$, $\mu_i$ is an 
integer.
\end{itemize} 

Let $T_i$ be the circle group with generator $v_i$, and $T'=\prod T_i$. 
Let $T'$ act on the product $M\times \C^r$, with moment map
\begin{equation}\psi_i(x,z)=\left\{\begin{array}{l@{\quad:\quad}l}
\l \Phi(x),v_i\r- \mu_i-|z_i|^2& i=1,\ldots, r\\
\l \Phi(x),v_i\r- \mu_i& i=r+1,\ldots,N.
\end{array}\right.\label{momentmap}
\end{equation}
For all $I\subset\{1,\ldots,N\}$, let $I'$ be the complementary set,
$R_I$ the subpolytope of $R$ defined by 
$\l\alpha,v_i\r=\mu_i$ for all $i\in I$, and $T_I=\prod_{i\in I}T_i$. 
Then 
\begin{equation}\psi^{-1}(0)\cong \bigcup_{I} \Phi^{-1}(\text{int}(R_I))
\times T_{I'}.\label{zerolevel}\end{equation}
Assumption (A) implies that for all $R_I\not =\emptyset$, 
$T_I$ is a finite cover of $T_{R_I}$, 
and by assumption (B), it acts locally freely on 
$ \psi^{-1}(\text{int}(R_I))$.
Hence $T'$ acts locally freely on $\psi^{-1}(0)$, and the cut space 
\begin{equation}M_R:=\psi^{-1}(0)/T'\end{equation}
is a symplectic orbifold. Moreover, we have a decomposition of $M_R$ 
into symplectic suborbifolds, 
\begin{equation}M_R=\bigcup_{R'<R}\,W_{R'}\label{union}\end{equation}
where 
\begin{equation}W_{R'}=\Phi^{-1}(\text{int}(R'))/T_{R'}.\end{equation}
\begin{remark}
\begin{enumerate}
\item  
For all faces $R_I\subset R$, there is a natural embedding 
$M_{R_I}\hra M_R$ as a symplectic suborbifold of codimension 
$2(\dim R-\dim R_I)$. 
The normal bundle $N_I$ of $\Phi^{-1}(\text{int}(R_I))/T_I$
in $M_R$ is given by the associated bundle 
$\Phi^{-1}(\text{int}(R_I))\times_{T_I}
\C^{|I|}$, with the standard action of $T_I$ on $\C^{|I|}$. 
Let us rewrite 
this in terms of $T_{R_I}$. The covering mapping $T_I\ra T_{R_I}$ 
gives rise to an isomorphism $\T_I\ra \T_{R_I}$. 
Let $\Lambda_I\subset \Lambda$ be 
the lattice generated by the $v_i,\,i\in I$, and $\alpha_i\in \T_{R_I}^*$ 
the preimage of the dual basis of $\Lambda_I^*$ under 
the isomorphism $\T_{R_I}^*\ra \T^*_I$. Then 
\begin{equation}N_I= \oplus_{i\in I}\,N_{-\alpha_i},
\label{normbdl}\end{equation}
where $N_{-\alpha_i}$ is the associated bundle 
\begin{equation}
\Phi^{-1}(\text{int}(R_I))\times_{T_{R_I}}\C\end{equation}
for the orbi-character $\exp(-2\pi i\l\alpha_i,\xi\r)$.
\item 
Given a symplectic action of a Lie group $G$ on $M$ which 
commutes with the $T$-action, one obtains a symplectic $G$-action 
on $M_R$. 
If the action is Hamiltonian, with moment map $J:M\ra\G^*$, 
the induced action on $M_R$ is Hamiltonian, and the moment map $J_R$
is obtained from the $T_I$-invariant restrictions 
$J|\Phi^{-1}(\text{int}(R_I))$. In particular, one has a Hamiltonian 
$T$-action on $M_R$, with moment polytope 
$\Delta_R=\Phi_R(M_R)=\Delta\cap R$. 
\item 
Cutting is local, in the sense that for all open subsets 
$U\subset \T^*$, there is a canonical symplectic isomorphism  
\begin{equation}\Phi_R^{-1}(U)\cong \Phi^{-1}(U)_R. \label{locality}
\end{equation}
In particular, one does not always need a global $T$-action to define 
the cut space. 
\item 
Given a $T\times G$-equivariant symplectic (resp. Hermitian) 
orbifold-vector bundle $E\ra M$, 
one obtains a $T\times G$-equivariant symplectic (resp. Hermitian) 
orbifold-vector bundle $E_R\ra M_R$, by letting 
$$  E_R=(E\boxtimes \C)|\Psi^{-1}(0)/T',$$ 
where we use the action of $T'$ on $E$ induced by the canonical 
map $T'\ra T$, and the trivial action on $\C$. 
More generally, one can let $T'$ act on $\C$ by a nontrivial 
character. (The gluing formulas proved in this paper won't depend 
on this choice.)
In particular, if $E=L$ is a pre-quantum line bundle, 
one will only obtain a pre-quantum bundle $L_R$ for $M_R$ if condition 
(Q) is satisfied and if one uses the character 
$\exp(-2\pi i\l\mu,\xi\r)$ 
for $T'$. In the sequel, we will always make this choice for 
pre-quantum line bundles, without mentioning this explicitly. 
 
For each face $R_I\subset R$, there is a natural identification
\begin{equation}E_R|\,\big(\Phi^{-1}(\text{int}(R_I))/T_I\big)=
\big(E|\Phi^{-1}(\text{int}(R_I))\big)/T_I.\label{isom}\end{equation}
If $E=L$ is a pre-quantum line bundle, this identification 
preserves the Hermitian structure and the connection. 
\item
The $T\times G$-equivariant Riemann-Roch numbers $\RR(M_R,E_R)=\chi_R$ 
depend only on the polytope $R$, not on the choice of the $v_i,\mu_i$. 
In fact,  
all $T$-fixed point contributions $\chi_{R,F}$ are independent of 
this choice.
\end{enumerate}
\end{remark}
\noindent{\bf Example:}\\ 
A compact
Hamiltonian $T$-space $M$ is called a {\em symplectic toric
orbifold} if $\dim M=2\dim T$ and the $T$-action is effective.
By a theorem of Delzant \cite{D88}, symplectic toric {\em 
manifolds} are completely classified by their moment polytopes 
$\Phi(M)=\Delta$. 
Lerman and Tolman \cite{LT94} have shown that in the orbifold 
case, any rational simple polytope $\Delta\subset\T^*$
can occur as a moment map image. To specify $M$,  one needs 
in addition a positive integer attached to each 
facet $\Delta_i$; this corresponds to  
choosing some lattice vector $v_i$ perpendicular to $\Delta_i$ . 
It was mentioned in \cite{L93} that every  symplectic toric
orbifold can be obtained by symplectic cutting: 
Let $M=\C P(1)^k$, with symplectic form 
on $\C P(1)$ equal to $l\in\R_{>0}$ times the Fubini-Study form. The 
moment map 
$$  \Phi_j(w_1,\ldots, w_k)=-l\,\f{|w_j|^2-1}{|w_j|^2+1}$$ 
($w_j\in\C\cup\{\infty\}$)
defines a Hamiltonian $(S^1)^k$ action on $M$, with moment map image 
the cube $P_l=\{\alpha\in (\R^k)^*|\,-l\le \alpha_j\le l\}$. 
Consider a 
{\em compact} polytope $R$ as above, and choose $l$ large enough 
such that $R\subset\text{int}(P_l)$. The corresponding cut 
space $M_R$ does not depend on $l$, and is the 
symplectic toric orbifold associated to $R$. 
If $l\in\N$, the $l$th power of the hyperplane bundle is a 
pre-quantum bundle for $(\C P(1),\,l\omega_{F.S.})$, and by taking 
exterior tensor products one obtains a pre-quantum bundle $L$ for 
$M$. The equivariant character $\chi=\RR(M,L)$ 
for the $T$-action on $M$ is given by 
$$ \chi(z)=\prod_{j=1}^k\Big(\sum_{\nu_j=-l}^l 
z_j^{\nu_j}\Big)=\sum_{\alpha\in P_l}\,z^\alpha,$$ 
for $z=(z_1,\ldots,z_k)\in T$. 
Hence, by applying the equivariant version of Corollary 
\ref{corcut} in stages, one recovers 
the well-known result that the equivariant Riemann-Roch number  
for $M_R$ is given by 
\begin{equation}\RR(M_R,L_R)(z)=\sum_{\alpha\in R}z^\alpha.
\label{Delzant}\end{equation}

\begin{remark}
Letting $D_R$ be the symplectic toric orbifold associated to $R$, 
and $M$ any Hamiltonian $T$-space, one can actually {\em define} 
$M_R$ to be the reduction at $0$ of $M\times D_R^-$ 
with respect to the diagonal action; here $D_R^-$ denotes $D_R$   
with the opposite symplectic structure. 
\end{remark}

\subsection{The abelian gluing formula}
Consider now a finite collection $\mathcal{R}=\{R\}$ of admissible 
polytopes, 
such that $\Delta\subset \bigcup_{R\in\mathcal{R}} R$, 
and such that for each polytope in $\mathcal{R}$, all 
faces are also in $\mathcal{R}$, and for all  
$R_1,R_2\in\mathcal{R}$, the intersection 
$R_1\cap R_2$ is a face of each. 
For each $R\in\mathcal{R}$, let $W_R=\Phi^{-1}(\text{int}(R))/T_R$. 
By (\ref{union}), each 
cut space $M_R$ is a disjoint union of symplectic orbifolds
\begin{equation}M_R=\bigcup_{R'\subset R}\,W_{R'}.\end{equation}
Although the gluing of the $W_{R'}$ depends on the choice of the 
$v_i$, the 
Riemann-Roch numbers $\RR(M_R,E_R)$ are independent of this choice.  

\begin{theorem} \label{abcut} (Gluing Formula)
The Riemann-Roch numbers satisfy the gluing rule 
\begin{equation}(-1)^{\dim\Delta}\,\RR(M,E)=\sum_{R\in\mathcal{R}}
(-1)^{\dim R}\,\RR(M_R,E_R).\end{equation}
If a compact Lie group $G$ acts on $E\ra M$ and this action 
commutes with the action of $T$, the gluing rule holds for the 
corresponding $G$-equivariant Riemann-Roch numbers.
\end{theorem}    
We can actually prove a slightly stronger, local result. Write 
each $\RR(M_R,E_R)=\chi_R$ as a sum over fixed point 
contributions $\chi_{R,F}$, where $F$ ranges over the connected 
components of $M_R^T$.  

Suppose $F\subset W_S$ is a connected component of 
a $T$-fixed point orbifold for some $S\in\mathcal{R}$. 
If $\dim S=\dim T$, $W_S$ is an open subset of 
$M$, so $F$ is also a fixed point 
orbifold for $M$, and clearly $\chi_F=\chi_{S,F}$.
Suppose, on the other hand, that $\dim S< \dim T$. Then $F$ 
is {\em not} a fixed point orbifold of $M$, and is a 
fixed point orbifold for $M_R$ if and only if 
$\Phi_S(F)\in R$.
\begin{theorem}\label{locglueform} 
(Local gluing formula)
Let $F\subset  \Phi^{-1}(\text{int}(S))/T_{S}$ be a fixed point orbifold,
where $S\in\mathcal{R}$ is a polytope with $\dim S<\dim T$. Then 
\begin{equation}\sum_{R\ni\mu}(-1)^{\dim R}\,\chi_{R,F}=0,\label{local}\end{equation}
where $\mu=\Phi_S(F)$.
\end{theorem}
Recall that for the $S^1$ case, we needed the identity 
$(1-z)^{-1}+(1-z^{-1})^{-1}=1$ in the proof of the gluing formula. 
To prove Theorem \ref{locglueform}, we need a 
somewhat more sophisticated version of this identity.

Consider a $k$-dimensional simplicial cone 
$C\subset \Lambda_\R^*:=\T^*$, and let $\alpha_j\in\Lambda^*
\otimes_\Z\Q$ be linearly independent generating vectors for $C$, 
with the property that the lattice $\Lambda_C^*$ generated by 
the $\alpha_j$ contains $\Lambda^*$. 
Let $\Lambda_C$ be the dual lattice, and $\Gamma_C=\Lambda/\Lambda_C$. 
For all $\gamma\in \Gamma_C$, let $c_j(\gamma)=
e^{2\pi i\l\alpha_j,\gamma\r}$, 
and define the following meromorphic function on the 
complexified torus $T^\C$:
\begin{equation}f_C(z)=\f{1}{\#\Gamma_C}\sum_{\gamma\in\Gamma_C}
\Big(\prod_{j=1}^k(1-c_j(\gamma)
z^{\alpha_j})\Big)^{-1},\label{fC}\end{equation}
where $z^\mu=\exp(2\pi i\l\mu,\xi\r)$ if $z=\exp(\xi)\in T^\C$. 
Notice that the individual summands on the right hand side are only
functions on the covering torus $(\Lambda_\R/\Lambda_C)^\C$, but 
the sum is $\Gamma_C$-invariant and therefore descends to $T^\C$. 
If $C$ is a lower dimensional simplicial polytope, $f_C$ is 
defined similarly, by considering the lattice
$(\R.C)\cap\Lambda^*$. 

\begin{lemma} \label{char}
On the set of all $z\in T^\C$ such that $|z^\mu|<1$ for 
all $\mu\in C-\{0\}$,   
\begin{equation}f_C(z)=\sum_{\mu\in C}\,z^\mu.\end{equation}
\end{lemma}
\begin{proof}
By Taylor's expansion, 
$$  f_C(z)=\f{1}{\#\Gamma_C}\sum_{\{\nu_j\ge 0\}}
\Big(\sum_{\gamma\in\Gamma_C}
e^{2\pi i\sum_j\nu_j\l\alpha_j,\gamma\r}\Big)\,
z\,^{\sum_j\nu_j\alpha_j}.$$ 
But the sum over $\Gamma_C$ equals $\#\Gamma_C$ if 
$\sum_j\nu_j\alpha_j\in\Lambda^*$, $0$ otherwise.\end{proof} 

Recall now that a simplicial fan in $\Lambda_\R^*$ is a 
collection $\lie{C}=\{C\}$ of simplicial cones, with the property 
that for each cone $C$ in $\lie{C}$, all faces are also in $\lie{C}$, 
and the intersection of any two cones in $\lie{C}$ is a face of each.
The fan is called complete if the $C$'s cover all of $\Lambda_\R$.
A fan $\lie{C}'$ is called a refinement of $\lie{C}$ if 
any cone in $\lie{C}$ is a union of cones in $\lie{C}'$.  
We define 
\begin{equation}f_\lie{C}(z)=\sum_{C\in\lie{C}}\,(-1)^{\dim C}\,f_C(z).\end{equation}
It follows from Lemma \ref{char} that $f_\lie{C}(z)$ 
is invariant under (simplicial) refinements. 

\begin{lemma}\label{poly}
If $\lie{C}$ is a complete simplicial fan,  
\begin{equation}\sum_{C\in\lie{C}}(-1)^{\dim C}\,f_C(z)=0.
\label{fanglu}\end{equation}
\end{lemma} 
\begin{proof}
Let $\beta_1,\ldots,\beta_k\in \Lambda^*$ be a lattice basis, 
and $\lie{C}_1$ the simplicial fan whose cones are generated 
by all subsets of $\pm \beta_1,\ldots,\pm\beta_k$ not 
containing both $\pm\beta_j$, for any $j$. Let $\lie{C}_2$ 
be the fan whose cones are all intersections of cones in 
$\lie{C},\lie{C}_1$, and let $\lie{C}_3$ be a 
{\em simplicial} refinement of $\lie{C}_2$. 
Since $\lie{C}_3$ refines both 
$\lie{C},\lie{C}_1$, it follows that $f_\lie{C}=f_{\lie{C}_3}
=f_{\lie{C}_1}$. But 
$$  f_{\lie{C}_1}(z)
=(-1)^k\prod_{j=1}^k\Big(\f{1}{1-z^{\beta_j}}+ 
\f{1}{1-z^{-\beta_j}}
-1\Big)=0,$$ 
q.e.d.
\end{proof}    

\noindent{\bf Proof of Theorem \ref{locglueform}.}
Without loss of generality, we may assume $\mu=0$.
We will also assume $G=\{e\}$ for simplicity, the general case is 
obtained by working 
with $T\times G$-equivariant characteristic classes. 

For each polytope $R\ni 0$, let $C=C_R$ be the cone 
$\R_{\ge 0}R$, and let $\lie{C}=\{C\}$ be the fan 
obtained in this way.  
Since the statement of the theorem depends only on local data 
in a neighborhood of $\Phi^{-1}(0)$, we can assume without loss 
of generality that $\mathcal{R}=\lie{C}$. 
Let us first consider the case $\dim{S}=0$, which implies that 
$S=\{0\}$ is a vertex of all $C\in\lie{C}$.  
This means in particular that $0$ is a regular value of $\Phi$, and 
that 
$F=M_{0}=\Phi^{-1}(0)/T$. By assumption (A),  
$\lie{C}$ is a complete simplicial fan in $\T^*$. 

We now have to investigate 
how the local contributions from the fixed 
point 
formula add up. Notice that the multiplicity of 
the fixed point manifold $F\subset M_C$ depends on 
$C$. In fact, the generic isotropy group of $F$ as a 
suborbifold of $M_C$ is a $\#\Gamma_C$-fold cover 
of the isotropy group of $F$ as identified with $M_S$, 
since this is the ``twist'' of the normal bundle 
$N_C$ of $F$ in $M_C$. It follows that the associated 
orbifold
$\ti{F}_C$ of $F$ considered as a suborbifold 
of $M_C$ is a $\#\Gamma_C$-fold cover of 
$\ti{F}:=\ti{F}_S$.

By (\ref{normbdl}), there is an 
explicit description of $N_C$: Suppose 
the cone $C$ is generated by the orbi-weights $\alpha_j$. 
For each orbi-weight $\alpha\in\Lambda^*\otimes_\Z\Q$, 
let $N_{\alpha}$ 
be the associated orbifold bundle 
$N_{\alpha}=\Phi^{-1}(0)\times_{T}\C$, where $T$ acts by 
the (orbi-) character $\exp(2\pi i\l\alpha,\xi\r)$.
Then $N_C=\oplus_{\alpha_j}\,N_{-\alpha_j}$. 

Consider now $E'=\Phi^{-1}(0)\times_T(\T\otimes\C)$ as a 
$T$-equivariant orbifold bundle over $F$, and let 
$F_\T({E'},\xi)$ be its equivariant curvature.
Let $\ti{E}'$ be the pullback of $E'$ to $\ti{F}$, 
and $A\in C^\infty(\text{Aut}(\ti{E}'))$ as  
in section \ref{RRO}. Notice that the equivariant 
curvature of $N_\alpha$ is 
$F_\T(N_\alpha,\xi)=\l\alpha,F_\T({E}',\xi)\r$.

By performing the summation over the fibers of 
$ \ti{F}_C\ra \ti{F}$, 
the fixed point contribution of $F\subset M_C$
becomes
$$  \chi_{C,F}(e^\xi)=
\int_{\ti{F}}
\f{1}{d_{\ti{F}}} \f{ \Td(\ti{F})
\,{\Ch}^{\ti{F}}_\T(\ti{E}|\ti{F},\xi)}
{D^{\ti{F}}(N_{\ti{F}})  }\,
f_C(A^{-1}\,e^{\f{i}{2\pi}F_\T(\ti{E}',\xi)}   ).
$$  
The local gluing formula now follows directly from 
(\ref{fanglu}).

If $S\not=\{0\}$, each $C$ contains $S$ as a subspace, and 
the collection of quotient cones $C/S$ defines a complete 
rational simplicial fan in $\T^*/S$.  
The normal bundle $N_C$ of $F$ in $M_{C}$ 
splits into the direct sum of its part in $M_{S}$ and its 
symplectic orthogonal, which we denote by $N_{C/S}$:
$$  N_C=N_{S}\oplus N_{C/S}.$$ 
Hence
$$  \chi_{C,F}(e^\xi)=
\int_{\ti{F}}
\f{1}{d_{\ti{F}}} \f{ \Td(\ti{F})
\,{\Ch}^{\ti{F}}_\T(\ti{E}|\ti{F},\xi)}
{D^{\ti{F}}(N_{\ti{F}})\,D^{\ti{F}}(N_S,\xi)}\,
f_{C/S}(A^{-1}\,e^{\f{i}{2\pi}F_\T(\ti{E}',\xi)}),
$$ 
where now $E'=\Phi^{-1}(\text{int}(S))\times_{T_S}(\T_S\otimes\C)$.
Again, the claim follows from (\ref{fanglu}).

\section{Nonabelian Cutting}

In this section, we will prove a generalization of Theorem \ref{abcut}  
to nonabelian groups. 
Let $G$ be a compact connected Lie group, with maximal torus $T$, and 
$G=K\,A$ 
its decomposition into its semisimple part, $K=(G,G)$, and its connected 
abelian part. 
By general properties of orbit type decompositions, the stabilizer 
group $G_\Sigma$ corresponding to a face $\Sigma$ is contained in 
$G_{\Sigma'}$ 
if and only if 
$\Sigma'$ is a face of $\Sigma$, with equality if and 
only if $\Sigma=\Sigma'$. 
Let $G_\Sigma=K_\Sigma\, A_\Sigma$ be the decomposition into 
semisimple and abelian part.
The Lie-algebras are given by $\K_\Sigma=
[\G_\Sigma,\G_\Sigma]$ and $\lie{a}_\Sigma=(\G_\Sigma)^{G_\Sigma}$, 
respectively. 
For all $\alpha\in \Sigma$, 
\begin{equation}T_\alpha(\Sigma)=(\G^*)^{G_\Sigma}\cap \T^*=
\lie{a}^*_\Sigma=[\G_\Sigma,\G_\Sigma]^0\cap \T^*.\end{equation}
Consider now a compact 
connected Hamiltonian $G$-orbifold $M$, with moment map 
$J:M\ra\G^*$. 
By a theorem of Kirwan \cite{K85}, 
$\Delta=J(M)\cap \T^*_+$ is a 
compact convex polytope. 
(To be precise, Kirwan's theorem only covers the manifold case, 
the extension to Hamiltonian orbifolds is proved in \cite{LMTW95}.)

Given a face $\Sigma$ of $\T^*_+$, let $U_\Sigma$ be the open 
set  
\begin{equation}U_\Sigma=\{\alpha\in\T^*_+|\,G_\alpha\subset G_\Sigma\}=
\bigcup_{\Sigma \subset\Sigma'}\text{int}(\Sigma'),\end{equation}
and define 
\begin{equation}Y_\Sigma=J^{-1}(G_\Sigma.U_\Sigma),\,\,M_\Sigma=J^{-1}(G.U_\Sigma)
=G\times_{G_\Sigma}Y_\Sigma.\end{equation}

Let $\pi_\Sigma:\T^*\ra \lie{a}^*_\Sigma$ denote the 
projection, and  $ q:\G^*\ra \T^*_+$ the mapping that sends $\alpha$ 
to the unique point of intersection 
of the coadjoint orbit $G.\alpha$ with $\T^*_+$. 
Notice that the restriction of $q_\Sigma:=\pi_\Sigma\circ q$ to 
$G.U_\Sigma$ is 
smooth.
\begin{theorem}\label{slicethm} 
(Symplectic cross section theorem.) $Y_\Sigma$ 
is a connected symplectic submanifold of $M$, and is a 
Hamiltonian $G_\Sigma$-space, with the restriction of $J$ serving as a 
moment map. The action of $A_\Sigma\subset G_\Sigma$ on $Y_\Sigma$ 
extends 
in a unique way to an action on $M_\Sigma$ which commutes with the 
$G$-action. Moreover, this action is Hamiltonian, with moment map 
$\Phi_\Sigma=q_\Sigma\circ J$. 
\end{theorem}
For a proof of the first part, see \cite{GS84}, p. 194. The second 
part is obvious since $\Phi_\Sigma=J$ on $Y_\Sigma$.

The idea to use the local $A_\Sigma$-actions for symplectic 
cutting of Hamiltonian $G$-spaces is due to Chris Woodward 
\cite{W95}. 
Consider a simple polytope $R\subset \T^*$ of the form 
(\ref{de}). 
Suppose that for all faces $S$ of $R$, and all $\Sigma$ 
such that $S\cap \Sigma\not=0$, the torus $T_{S}$ is a  
subtorus of $A_\Sigma$. 
By taking perpendiculars, the condition 
\begin{equation}S\cap \Sigma\not=\emptyset \Longrightarrow 
T_{S}\subset {A}_\Sigma\end{equation}
for all $S,\Sigma$ is equivalent to the following assumption:
\begin{itemize}
\item[(C)] 
For all faces $S$ of $R$ meeting a face $\Sigma$ of $\T^*_+$, 
the tangent space to $S$ contains the space perpendicular 
to $\Sigma$ (i.e. the space $\K_\Sigma^*\cap \T^*$). 
\end{itemize}
For all $\Sigma$, choose a neighborhood $V_\Sigma\subset U_\Sigma$
such that $V_\Sigma\cap R=V_\Sigma\cap\pi_\Sigma^{-1}(R_\Sigma)$
where $R_\Sigma=R\cap{\lie{a}_\Sigma}^*$, and define the cut space 
$M_R$ by gluing the cut spaces with respect to the 
local $A_\Sigma$-actions, $(J^{-1}(G.V_\Sigma))_{R_\Sigma}$. 
$M_R$ is a symplectic orbifold if for all $\Sigma$, $R_\Sigma$ is 
admissible for $J^{-1}(G.V_\Sigma)$, in other words if 
the following condition is satisfied:
\begin{itemize}
\item[(B')] 
For all faces $S$ of $R$ and all 
$x\in J^{-1}(S\cap\T^*_+)$, $\G_x\cap \T_{S}=\{0\}$.
\end{itemize}
To summarize, in the nonabelian case we will call a polytope 
$R\subset\T^*$ {\em admissible} if it satisfies conditions 
(A), (B') and (C), plus the extra condition (Q) if we 
are in the pre-quantization setting.  
Each admissible polytope $R$ defines a cut space $M_R$. 
The discussion from the abelian case goes through with 
no essential changes: 
We have a decomposition into symplectic sub\-or\-bi\-folds,   
\begin{equation}M_R=\bigcup_{S\subset R}\, 
(q\circ J)^{-1}(\text{int}(S))/T_{S}, \end{equation}
and for each face $S\subset R$, there is a canonical embedding 
$M_{S}\ra M_R$ as a symplectic suborbifold of 
codimension $2(\dim R-\dim S)$.
The induced action of $G$ on $M_R$ is Hamiltonian, 
and the corresponding moment polytope is 
\begin{equation}\Delta_R=J_R(M_R)\cap\T^*_+=\Delta\cap R.\end{equation}
If $E\ra M$ is a $G$-equivariant complex vector bundle, 
one obtains a cut-bundle $E_R\ra M_R$ with a canonically  
induced $G$-action.\\ 

\begin{remark}
A more geometric description of condition (B') can be obtained 
as follows. For all subalgebras $\H\subset\G$, let $M_{(\H)}\subset 
M$ be the set of points $x$ with 
infinitesimal stabilizer $\G_x$ conjugate to $\H$. It is 
well-known that only finitely many conjugacy classes $(\H)$ 
occur as stabilizers, and that each $ M_{(\H)}$ has a 
finite number of connected components. 
For each representative 
$\H$ for $(\H)$, let $M_\H$ be the set of all 
$x$ with $\G_x=\H$; then $M_{(\H)}=G.M_\H$.
By equivariance of the moment map, 
\begin{equation}J(M_\H)\subset (\G^*)^\H=\{\alpha\in\G^*|\,\text{ad}^*(\xi)
\alpha=0\}.\label{dfl}\end{equation} 
Let $\lie{z}=\G^\H$ be the centralizer of $\H$ in $\G$, 
and $Z=\exp(\lie{z})$. Then $Z\subset G$ acts on 
$M_\H$ in a Hamiltonian fashion, and by (\ref{dfl}) 
the restriction $J|M_\H$ serves as a moment map. 
Since every moment map $J$ has the property  
\begin{equation}\text{im}(T_x\,J)=\G_x^0 \mbox{ for all }x \in M,\end{equation}
this shows that the image under $J$ of each connected component 
of $M_\H$ is an open subset of an affine space of the form 
$\alpha+(\H^0)^\H$, where $\alpha\in \lie{z}^*$. 
By a suitable choice of $\H$, we can assume 
that $\T_1=\lie{z}\cap\T$ is a Cartan subalgebra 
of $\lie{z}$. Then 
$J(M_{\H})\cap\T^*=J(M_{\H})\cap\T_1^*$, and therefore 
$$  J(M_{(\H)})=G.J(M_{\H})=G.Z.(J(M_{\H})\cap\T_1^*)
=G.(J(M_{\H})\cap\T^*).$$  
It follows that the image under $\Phi$ of  
each connected component of $M_{(\H)}$ 
is an open connected subset of  
$W.(\alpha+(\H^0\cap \T^*))\cap \T^*_+$ for some 
$\alpha\in\lie{z}^*\cap \T^*$.
Using this result, condition (B') can be formulated as follows: 
\begin{itemize}
\item[(B)] 
The faces of $R$ intersect 
all sets $q\circ J(M_{(\H)})$ transversally.
\end{itemize}
Notice that this condition is generically fulfilled. 
\end{remark}
\begin{theorem} (Nonabelian gluing formula)
Suppose that $\mathcal{R}=\{R\}$ is a finite collection of 
admissible polytopes 
such that the $R$'s cover $\Delta$, the faces of each 
$R\in \mathcal{R}$ are also in $\mathcal{R}$, 
and for all $R_1, R_2\in\mathcal{R}$, 
their intersection is a face of each.
For each $G$-equivariant vector bundle $E\ra M$, 
one has the following gluing formula 
for $G$-equivariant Riemann-Roch numbers: 
\begin{equation}(-1)^{\dim\Delta}\,\RR(M,E)=\sum_{R\in\mathcal{R}}
(-1)^{\dim R}\,\RR(M_R,E_R).\end{equation}
\end{theorem}
\begin{proof}
This follows again from a local result:
Suppose that $F$ is a $T$-fixed point manifold for some $M_R$, 
and  $\alpha=q\circ J_R(F)$. If $\dim R=\dim T$ and 
$\alpha\in\text{int}(R)$, then $F$ is a $T$-fixed point manifold 
for $M$, and $\chi_F=\chi_{R,F}$. 
Otherwise, let $S\in\mathcal{R}$ be the unique polytope with 
$\alpha\in\text{int}(S)$. Since $M_S\subset M_R$ if $S\subset R$, 
$F$ is $T$-fixed for all $M_R$ with $\alpha\in R$. We have to show that  
\begin{equation}\sum_{R\ni\alpha}(-1)^{\dim R}\,\chi_{R,F}(e^\xi)=0.\label{eq23}\end{equation}
The normal bundle $\nu_{R,F}$ of $F$ in $M_R$ is a direct sum of 
the normal bundle of $F$ in $M_S$, and 
the restriction to $F$ of the normal bundle of 
$M_S$ in $M_R$.  Hence (\ref{eq23}) follows from the abelian result, 
Theorem \ref{locglueform}.
\end{proof}

We will now explicitly describe a decomposition 
$\mathcal{R}$ of 
$\T^*$ into admissible polytopes, which we will use in the following 
section. 
Assume first that $G$ is semisimple.
Let $\lie{S}_+=\{\beta_1,\ldots,\beta_k\}
\subset\lie{R}_+$ be the set of simple positive roots,
and $\alpha_1,\ldots,\alpha_k\in \Lambda^*\otimes_\Z\Q$ 
generating vectors for the edges of $\T^*_+$, such that $\beta_i$ 
is perpendicular to the facet spanned by the $\alpha_j$, $j\not=i$. 
Let $\lie{C}$ be the complete simplicial fan in $\T^*_+$ generated 
by $\alpha_1,\ldots,\alpha_k,-\beta_1,\ldots,-\beta_k$. 
The cones in this fan are generated by all subsets which do 
not contain both $\alpha_i$ and $-\beta_i$, for any $i$. 
For generic choices 
$\gamma\in\text{int}(\T^*_+)\cap\Lambda^*\otimes_\Z\Q$,
the polytopes $R_C=\gamma+C$ are admissible. Notice that for
the polytope $R_0$ that contains $0$, the intersection 
$R_0\cap\T^*_+$ is equal to the intersection of the convex hull 
of $W.\gamma$ with $\T^*_+$.     
The following picture shows the decomposition for the Weyl 
chamber of $G=SU(3)$, the intersection $R_0\cap\T^*_+$ is shaded. 

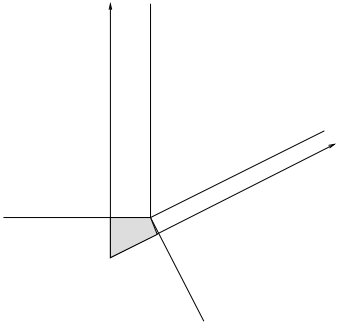
\begin{center}
\setlength{\unitlength}{0.0035in}
\begin{picture}(495,490)(0,-10)
\path(162.000,467.000)(160.000,475.000)(158.000,467.000)
\path(160,475)(160,95)(495,265)
\path(488.771,259.596)(495.000,265.000)(486.961,263.163)
\path(0,155)(220,155)(300,0)
\path(220,475)(220,155)(480,285)
\texture{c0c0c0c0 0 0 0 0 0 0 0 
	c0c0c0c0 0 0 0 0 0 0 0 
	c0c0c0c0 0 0 0 0 0 0 0 
	c0c0c0c0 0 0 0 0 0 0 0 }
\shade\path(160,155)(220,155)(230,130)
	(160,95)(160,155)
\end{picture}

\end{center}

If $G$ is a general compact connected group, 
let $G=K\,A$ be the decomposition into semisimple and abelian 
part. We can apply the above construction to the semisimple part, 
to obtain a decomposition $\mathcal{R}^{(1)}$ of $\K^*\cap\T^*$ into 
polytopes. 
Now  take any decomposition $\mathcal{R}^{(2)}$
of $\lie{a}^*$ into rational polytopes satisfying assumption (A),
and let $\mathcal{R}$ be the set of all products $R^{(1)}\times R^{(2)}$, 
$R^{(i)}\in\mathcal{R}^{(i)}$.  
Again, these polytopes will be admissible for generic choices 
$\gamma\in \K^*\cap\text{int}(\T^*_+)$. 

The important point about this decomposition is that for all polytopes 
$R=R^{(1)}\times R^{(2)}$ with $0\not\in R^{(1)}$, there exists
a face  $\Sigma\not=\a^*$ of $\T^*_+$ such that the intersection 
$R\cap \T^*_+$ is contained in $ U_\Sig\cap\T^*_+$. This means  
that the cut space $M_R$ admits a global cross-section 
$Y_{R,\Sig}$:
$$ M_R=G\times_{G_\Sigma}Y_{R,\Sig}.$$   
In this situation, Theorem \ref{slice} shows that 
the computation of $\RR(M_R,E_R)$ reduces 
to that of $\RR(Y_{R,\Sig},E_R| Y_{R,\Sig})$.

\section{Quantization}

Up to this point, we have been dealing with arbitrary $G$-equivariant
vector bundles $E$. Let us now focus on the special case that $E=L\ra
M$ is a pre-quantum line bundle, as explained in the introduction. Let
$\chi=\RR(M,L)\in R(G)$ denote the equivariant Riemann-Roch number,
and $N:\,\Lambda^*_+\ra\Z$ the multiplicity function. In this section,
we will prove the Guillemin-Sternberg conjecture, Theorem
\ref{GSC}. We will use that Theorem \ref{GSC} is already proved in the 
abelian
case (see the remarks at the end of section \ref{SymplecticSurgery}).
Using induction on $\dim G$, we can also assume that Theorem \ref{GSC}
(hence also Corollary \ref{SupportMultiplicity}) holds for all proper
subgroups of $G$. 
Choose a decomposition $\mathcal{R}=\mathcal{R}^{(1)}\times \mathcal{R}^{(2)}$ 
 as described at the end of the previous section, in such a way that $0$ is 
contained in the interior of a unique polytope $R_0=R^{(1)}_0\times R^{(2)}_0$ 
We claim that for any polytope $R=R^{(1)}\times R^{(2)}$ with 
$0\not\in R$, we have $\RR(M_R,L_R)^G=0$. Indeed, if $0\not\in R^{(2)}$ 
this follows from the abelian result since already $\RR(M_R,L_R)^A=0$. 
If $0\not\in R^{(1)}$, we have a global symplectic cross-section 
$Y_{R,\Sig}$, and the result follws from Theorem \ref{slice} because 
$\RR(Y_{R,\Sig},E_R| Y_{R,\Sig})^{G_\Sig}=0$ since 
$\Phi(Y_{R,\Sig})\not\ni 0$. Since we may choose $R_0\cap\T^*_+$ 
arbitrarily small, this proves the follwing excision property for 
$\RR(M,L)^G$:  
\begin{proposition}
If $R\subset \T^*$ is any admissible polytope containing $0$, 
$\RR(M,L)^G=\RR(M_R,L_R)^G$. 
\end{proposition} 
In particular, $\RR(M,L)^G$ depends only upon local data near $J^{-1}(0)$.
Let us now suppose that $0\in J(M)$ is a regular value of $J$, and denote 
$P=J^{-1}(0)$. Then $\pi:\,P\ra M_{red}$ is an orbifold-principal 
$G$-bundle. 
We will use the normal form theorem of Gotay to describe a neighborhood 
of $P$ in $M$.
Choose a connection $\theta\in\A^1(P,\G)$ on $P$, 
and equip the product $P\times \G^*$ with the closed two-form
\begin{equation}\sigma=\pi^*\omega\,+\d\,\l\text{pr}_2,\theta\r,\end{equation}
where $\text{pr}_2:\,P\times \G^*\ra\G^*$ 
denotes projection to the second factor. 
The diagonal action of $G$ makes $P\times \G^*$ into a Hamiltonian
$G$-space, with moment map equal to $\text{pr}_2$. 
\begin{theorem} (Normal form theorem \cite{G82}.) 
There exists a $G$-equivariant 
symplectomorphism from a neighborhood of $P$ in $P\times\G^*$ 
to a neighborhood of $P$ in $M$.
\end{theorem}
A well-known consequence of the normal form theorem is 
the following description of reduced spaces, for 
$\alpha$ close to $0$. Let $F^\theta=\d\theta\in \A^2(P,\G)$ 
be the curvature of $\theta$. Let $\mathcal{O}_\alpha$ be 
the coadjoint orbit through $\alpha$, equipped with its 
canonical symplectic structure $\sigma_\alpha$. 
The $G$-action on $\mathcal{O}_\alpha$ 
is Hamiltonian, with moment map the embedding 
$J_\alpha:\,\mathcal{O}_\alpha\hra \G^*$.  
By the shifting trick, 
the reduced space $M_\alpha=J^{-1}(\alpha)/G_\alpha$ 
is symplectomorphic to the reduced space at zero of 
$M\times \mathcal{O}_\alpha^-$, where the superscript ``-'' indicates 
that one takes the opposite symplectic form.
By doing this 
calculation in the canonical model, one finds: 
\begin{corollary}\label{redsp}
There is a neighborhood $U\ni 0$ in the set of regular values of 
$J$, such that for all $\alpha\in U$, the reduced space 
$M_\alpha$ is symplectomorphic to the symplectic fiber bundle 
$$  M_\alpha=P\times_G\,\mathcal{O}_\alpha\stackrel{\phi_\alpha}{\lra} 
M_0,$$ 
with symplectic form given by the minimal coupling recipe 
of Sternberg \cite{S77}
\begin{equation}\pi_\alpha^*\omega_\alpha=\pi_0^*\omega_0\,+\d\l J_\alpha,\theta\r
-\sigma_\alpha.
\label{mincoup}\end{equation}
Here $\pi_\alpha:\,P\times_G\,\mathcal{O}_\alpha\ra M_\alpha$ 
denotes the projection.
\end{corollary}
Similarly, one can express the associated orbifold
$\ti{M}_\alpha$ as a symplectic fiber bundle 
\begin{equation}
\ti{M}_\alpha\cong
\hat{P}\times_G \mathcal{O}_\alpha\stackrel{\ti{\phi}_\alpha}{\lra}
\ti{M}_0,\end{equation}
where $\hat{P}$ was defined  in (\ref{hatP}). Notice that 
$N_{\ti{M}_\alpha}=\ti{\phi}_\alpha^* N_{\ti{M}_0}$, and hence
\begin{equation}D^{\ti{M}_\alpha}(N_{\ti{M}_\alpha})
=\ti{\phi}_\alpha^* D^{\ti{M}_0}(N_{\ti{M}_0})\end{equation}

\noindent{\bf Proof of Theorem \ref{GSC}:} 
Choose $R_0=R_0^{(1)}\times R_0^{(2)}$ as above, in such a 
way that $R_0\cap\T^*_+$ is contained in the neighborhood $U$
from Corollary \ref{redsp}.
The moment polytope for $M_{R_0}$ is then simply 
$\Delta_{R_0}=\T^*_+\cap R_0$, and the set of regular values 
of $J_{R_0}$ is $\text{int}(R_0)\cap\T^*_+$. 
Notice that $0$ is a regular value 
for the action of $A\subset G$ on $M_{R_0}$, since $\K^*\oplus\{0\}$ is 
transversal to all faces of $R_0$. From the result for the abelian 
case, we have 
$$  \RR(M,L)^G=\RR(M_{R_0},L_{R_0})^G=\RR((M_{R_0})_A,(L_{R_0})_A)^K,$$ 
where the subscript $A$ denotes the reduced space with respect to 
the $A$-moment map.
Using reduction in stages, it is therefore sufficient to prove Theorem 
\ref{GSC} for the semisimple case. 

Let us assume for the rest of this proof that $G$ is semisimple.
Let $\{\beta^\sharp|\,\beta\in\lie{S}_+\}\subset 
\Lambda\otimes_\Z\Q$ be the dual basis to $\lie{S}_+$.  
By definition, $R_0$ is given by inequalities 
\begin{equation}\l\alpha-\gamma,\beta^\sharp\r\le 0
\,\,\mbox{ for all }\beta\in\lie{S}_+.\end{equation}
Choose $k\in\N$ such that for all $\beta\in\lie{S}_+$, 
$v_\beta:=-k\beta^\sharp\in\Lambda$ and $\mu_\beta:=-k\l\gamma,\beta^
\sharp\r
\in\N$. Thus $R_0$ is given by 
\begin{equation}\l\alpha,v_\beta\r\ge \mu_\beta\,\,\mbox{ for all }
\beta\in\lie{S}_+.\end{equation}
To compute $N(0)$, we may replace $L\ra M$ by the cut bundle 
$L_{M_{R_0}}\ra  M_{R_0}$, which we continue to denote by 
$L\ra M$.   
The components of the $T$-fixed point set for $M$ are then simply the 
preimages 
\begin{equation}J^{-1}(w.\gamma)\cong M_{w.\gamma}\cong M_\gamma\cong P/T,\end{equation}
for $w\in W$. 
We can write $\chi(e^\xi)$ as a sum 
over fixed point contributions $\chi_{\gamma}$ of $M_\gamma$ considered 
as a fixed point manifold in the symplectic cross-section 
$J^{-1}(\text{int}(\T^*_+))$:
\begin{equation}\chi(e^\xi)=\f{ 
\sum_{w\in W}\det(w)\,e^{2\pi i \l w(\delta)-\delta,\xi\r} 
\chi_{\gamma}(w^{-1}(\xi)) }
{\prod_{\beta\in\lie{R}_+}(1-e^{-2\pi i\l\beta,\xi \r})}\end{equation}
The fixed point formula for $\chi_{\gamma}$ involves the 
normal bundle $\nu_\gamma$ of $M_\gamma$ in  
$J^{-1}(\text{int}(\T^*_+))$. 
By (\ref{normbdl}),
$\nu_\gamma=\oplus N_{-\alpha_\beta}$, where the orbi-weights 
$\alpha_\beta:=-\f{1}{k}\beta$ form the dual basis to 
$v_\beta$.
The equivariant Chern class of $N_{-\alpha_\beta}$ is simply 
$c_\beta(\xi)=-\l \alpha_\beta,F^\theta+2\pi i\xi\r$. 
We have
\begin{equation}\chi_{\gamma}(w^{-1}(\xi))
=\int_{{\ti{M}_\gamma}}\f{1}{d_{{\ti{M}_\gamma}}}\,
\f{   \Td({\ti{M}_\gamma})\,\Ch_\T^{{\ti{M}_\gamma}}
(\ti{L}_{w(\gamma)},\xi)      }
{  D^{{\ti{M}_\gamma}}(N_{{\ti{M}_\gamma}})       }
\,\,
\prod_{\beta\in\lie{S}_+}\,D_\T^{\ti{M}_\gamma}
(N_{-w(\alpha_\beta)},\xi)^{-1}.\label{eq75}
\end{equation}
Let us consider (\ref{eq75}) as an equality  of meromorphic functions 
on $\T\otimes\C$. 
For all $w(\alpha_\beta)$, we can expand 
\begin{equation}
D_\T^{\ti{M}_\gamma}(N_{-w(\alpha_\beta)},\xi)^{-1}=
(1-e^{2\pi i\l{w(\alpha_\beta)},\xi\r}\,\Ch^{\ti{M}_\gamma}
(N_{{w(\alpha_\beta)}}))^{-1}
\end{equation}
into a geometric series with respect to $e^{2\pi
i\l{w(\alpha_\beta)},\xi\r}$.  Of course, this will only converge if
$\l{w(\alpha_\beta)},\,\lie{Im}(\xi)\r>0$.  Since we want to expand
all factors in (\ref{eq75}) simultanously we polarize the weights: 
Let
\begin{equation}
l_\beta^w=
\left\{\begin{array}{l@{\quad:\quad}l}
0& \l w(\alpha_\beta),\eta\r > 0\\
1& \l w(\alpha_\beta),\eta\r < 0
\end{array}
\right. 
\text{ for all }\eta\in\text{int}(\T_+),
\end{equation}
where $\T_+$ is the positive Weyl chamber in $\T$,
and write $\alpha_\beta^w=(-1)^{l_\beta^w}\,w(\alpha_\beta)$,
and
$\epsilon_w=(-1)^{\sum l_\beta^w}$.
Then we may rewrite the denominator of the second term 
in (\ref{eq75}), for $\lie{Im}(\xi)\in\text{int}(\T_+)$, as 
\begin{eqnarray*} 
D_\T^{\ti{M}_\gamma}(N_{-w(\alpha_\beta)},\xi)^{-1}
&=&(-1)^{l_\beta^w}\,D_\T^{\ti{M}_\gamma}(N_{-\alpha_\beta^w},\xi)^{-1}
\,\Ch_\T^{\ti{M}_\gamma}(N_{l_\beta^w\,\alpha_\beta^w},\xi)
\\&=&(-1)^{l_\beta^w}\sum_{l_\beta\ge 0} 
\Ch_\T^{\ti{M}_\gamma}(N_{(l_\beta+l_\beta^w)\,
\alpha_\beta^w},\xi).
\end{eqnarray*}
Notice also that $L_\gamma\cong\phi^*\,L_0\otimes N_\gamma$ 
by the pre-quantum condition and (\ref{mincoup}), thus 
$\Ch_\T^{\ti{M}_\gamma}(\ti{L}_\gamma,\xi)=
\ti{\phi}_\gamma^*\,\Ch^{\ti{M}_0}(\ti{L}_0)\,\,
\Ch_\T^{\ti{M}_\gamma}
(N_\gamma,\xi)$. We hence obtain the formula 
$$
\chi(e^\xi)\,\prod_{\beta\in\lie{R}_+}
(1-e^{-2\pi i\l\beta,\xi \r})        
=
\sum_{w\in W} \sum_{l_\beta\ge 0} 
\epsilon_w \,\det(w)\times$$
$$
\times 
e^{2\pi i\l w(\delta)-\delta,\xi\r}
\int_{\ti{M}_\gamma}
\f{1}{d_{\ti{M}_\gamma}}
\f{
\Td(\ti{M}_\gamma)\,
\ti{\phi}_\gamma^*\,\Ch^{\ti{M}_0}(\ti{L}_0)
}
{  
\ti{\phi}_\gamma^* D^{\ti{M}_0}(N_{\ti{M}_0})
}\,\,
\Ch_\T^{\ti{M}_\gamma}(N_{w(\gamma)+\sum(l_\beta+l_\beta^w)\,
\alpha_\beta^w}
,\xi).$$
On the other hand, Weyl's character formula
\begin{equation}\chi(e^\xi)\,  \prod_{\beta\in\lie{R}_+}
(1-e^{-2\pi i\l \beta,\xi\r})=
\sum_{\mu\in\Lambda^*_+} N(\mu){\sum_{w\in W}\,e^{2\pi i\l 
w(\delta+\mu)-\delta,\xi\r}}.
\label{Weyl}\end{equation}
shows that $N(\mu)$ is the coefficient of 
$e^{2\pi i\l\mu,\xi\r}$ in $\chi(e^\xi)\,\prod(1-e^{-2\pi i\l 
\beta,\xi\r})$
( because, for all $\nu\in\Lambda^*_+$, 
$w(\delta+\nu)-\delta\in\Lambda^*_+$ if and only $w=1$). 
To find $N(0)$, we thus have to solve the equation
\begin{equation} w(\delta+\gamma)-\delta =-\sum (l_\beta+l_\beta^w)\alpha_\beta^w,
\,\,l_\beta\ge 0.
\label{eq60}\end{equation}
Let us apply $w^{-1}$ to both sides, and take 
the scalar product with $\beta^\sharp$, using that 
$\alpha_\beta=-\f{1}{k}\beta$: 
\begin{equation} \l\delta+\gamma-w^{-1}(\delta),\beta^\sharp\r=
(-1)^{l_\beta^w}(l_\beta+l_\beta^w)\,k^{-1}, \,\, l_\beta\ge 0
\label{eq61}\end{equation}
The left hand side of (\ref{eq61}) is strictly positive, since 
$\delta-w^{-1}(\delta)$ is a sum of positive roots. 
If $w\not=1$, there is no solution of (\ref{eq60}), 
because at least one $l_\beta^w=1$, which  makes the right hand side 
negative. If $w=1$, we have $l_\beta^w=0$ for all $\beta\in\lie{S}_+$, 
and therefore $l_\beta=k \l \gamma,\beta^\sharp\r$. 
But then 
\begin{equation}\sum_{\beta\in\lie{S}_+}\,l_\beta\,\alpha_\beta
=
\sum_{\beta\in\lie{S}_+} \l\gamma,\beta^\sharp\r\, \beta
=\gamma.\end{equation}
We have thus shown: 
\begin{equation}
N(0)=\int_{\ti{M}_\gamma}\f{1}{d_{\ti{M}_\gamma}}
\f{\Td(\ti{M}_\gamma)\,
\ti{\phi}_\gamma^*
{\Ch}^{\ti{M}_0}
(\,\ti{L}_0)}
{  
\ti{\phi}_\gamma^* D^{\ti{M}_0}(N_{\ti{M}_0})}.
\end{equation}

Let us finally integrate over the fibers of $\ti{\phi}_\gamma$, 
using that 
$\int_{\mathcal{O}}\Td(\mathcal{O})=1$ for every coadjoint orbit 
$\mathcal{O}$: 
$$
N(0)=\int_{\ti{M}_0} \f{1}{d_{\ti{M}_0}}
\f{\Td(\ti{M}_0)\,{\Ch}^{\ti{M}_0}
(\ti{L}_0)}{D^{\ti{M}_0}(N_{\ti{M}_0})}=\RR(M_0,L_0),
$$
q.e.d.

\section{Appendix: A Short Proof for $G=SU(2)$}   
In this section, we will give a short proof of Theorem \ref{GSC} for 
the case $G=SU(2)$, modeled after the proof for $G=S^1$ in 
\cite{DGMW95}. The main idea will be to construct a Hamiltonian 
$S^1$-space, which has the same multiplicities and the same reduced
spaces.
Let $(M,\omega)$ be a quantizable Hamiltonian $G$-space, 
and suppose that $0$ is a regular value of the 
moment map $J$. 

For simplicity, we will assume that $M$ is a manifold, and that the 
action on 
$J^{-1}(0)$ is free, although the proof is easily adaptable to the 
orbifold case.
Let $T=S^1$ be the maximal torus of $SU(2)$. The dominant weights of $G$ 
are 
labeled by nonnegative integers, $\Lambda_+=\Z_{\ge 0}$, and the 
positive 
root is equal to $2$. By Weyl's character formula, 
\begin{equation}
\chi(e^{i\phi})=\sum_{\mu\in\Z_{\ge 0}}\,N(\mu)\,
\Big(\f{e^{i\mu\phi}}{1-e^{-2i\phi}}+\f{e^{-i\mu\phi}}{1-e^{2i\phi}}\Big).
\label{W}\end{equation} 
This shows that
the restriction of $\chi$ to $S^1$ extends to a rational function on the 
Riemann 
sphere, $\C\cup\{\infty\}$, and 
\begin{equation}N(0)=\text{res}_{z=\infty}\Big(\f{(1-z^{-2})\chi(z)}{z^{-1}}\Big).\end{equation}
On the other hand, we can use the equivariant index theorem 
(\ref{locc}) to express $\chi|S^1$ as a sum over fixed point 
contributions, $\chi|S^1(z)=\sum_F \chi_F(z)$. 
By analytic continuation, this becomes an equality of rational functions 
on $\C\cup\{\infty\}$.

Since $0$ is a regular value of $J$, all $S^1$-fixed
point components $F$ will have $J_F\not=0$. By the symplectic 
cross section 
theorem, $Y_+=J^{-1}(\R_{>0})$ is a symplectic submanifold of $M$, 
equipped with a Hamiltonian action of $T=S^1$ whose moment map is simply 
the restriction of $J$. The restriction of $L$ to $Y_+$ is a pre-quantum 
line 
bundle. 
For all $J_F>0$, we can view $F\subset Y_+$ as a fixed point 
manifold for the $S^1$ action on $Y_+$. Let 
$\chi_{+,F}(e^{i\phi})$ be
the corresponding fixed point contribution.
Then 
$$  \chi_F(z)=\f{\chi_{+,F}(z)}{1-z^{-2}}.$$ 
Note that the Weyl group $W=\{e,g\}\cong\Z_2$ 
of $G$ acts effectively on the set of 
connected components of $M^T$, and $\chi_{g.F}(z)=\chi_F(z^{-1})$.  
Hence 
\begin{equation}\chi(z)=\sum_{J_F>0}\Big(\f{\chi_{+,F}(z)}{1-z^{-2}}+
\f{\chi_{+,F}(z^{-1})}{1-z^{2}}\Big)
\end{equation}
or 
\begin{equation}(1-z^{-2})\,\chi(z)=\sum_{J_F>0}\big(\chi_{+,F}(z)-
z^{-2}\chi_{+,F}(z^{-1})
\big).\end{equation}
Since $\chi_{+,F}(z^{-1})=O(z^{-J_F})$ for $z\ra\infty$, it follows that only 
the first term will contribute to the residue at $z=\infty$, and we 
obtain the formula 
\begin{equation}N(0)=\sum_{J_F>0}\text{res}_{z=\infty}\,\f{\chi_{+,F}(z)}
{z^{-1}}.\end{equation}
We will now use the following trick, which is due to Eugene Lerman and may
be regarded as the $SU(2)$ version of symplectic cutting. 
Let $\C^2$ be equipped 
with its standard symplectic structure and the standard action of 
$U(2)$. Let $\phi:\,M\times (\C^2)^-\ra \lie{su}(2)^*$
be the moment map for the diagonal $SU(2)\subset U(2)$ 
action, and define $M_+$ to be 
the reduced space, 
$M_+=\phi^{-1}(0)/SU(2)$. It is easy to check the following properties 
of $M_+$:
\begin{enumerate}
\item 
As a set, $M_{+}$ is equal to the disjoint union $M_0\cup\, Y_+$, and the 
embeddings of $M_0$ and $Y_+$ are symplectic. 
\item 
The normal bundle of $M_0$ in $M_+$ is isomorphic to the associated 
bundle  
$$ \nu=J^{-1}(0)\times_{SU(2)}\C^2.$$
\item 
From the action of the center $U(1)\subset U(2)$ on $\C^2$ we 
get an induced Hamiltonian 
$U(1)$ action on $M_+$, which fixes $M_0$ and is equal to the action of
the maximal torus
$T\subset SU(2)$ on $Y_+$. The weights for the $U(1)$ action on $\nu$ 
are $(-1,-1)$. 
\item
Let $L_{\C^2}$ be the trivial line bundle over $\C^2$, with fiber metric 
$\exp(-\pi\,|z|^2)$. Then $L_+=(L\boxtimes L_{\C^2}^*|\phi^{-1}(0))/SU(2)$ 
is a pre-quantum line bundle for $M_+$. There is a natural lift of 
the $U(1)$ 
action on $M_+$ to $L_+$. 
\end{enumerate} 
Let $\chi_+=\RR(M_+,L_+)$ be the equivariant Riemann-Roch number of 
$M_+$ with 
respect to this $U(1)$ action, and $N_+(\mu)$ the multiplicity function.
By the fixed point formula for $\chi_+$, we have 
\begin{equation}\chi_+(z)=\sum_{J_F>0}\,\chi_{+,F}(z)
+\int_{M_0}\f{\Td(M_0)\,
\Ch(L_0)}
{\det(1-z\,e^{\f{i}{2\pi}F(\nu)})}.\label{7.1}\end{equation}
The second term is $O(z^{-2})$  for $z\ra \infty$ since $\dim_\C\nu=2$. 
Hence
\begin{equation}N_+(0)=\text{res}_{z=\infty}\f{\chi_+(z)}{z^{-1}}
=\sum_{J_F>0} \text{res}_{z=\infty}\f{\chi_{+,F}(z)}{z^{-1}}=
N(0).\end{equation}
On the other hand, we see from (\ref{7.1}) that $\chi_+(z)$ is 
holomorphic 
for $z\ra 0$, and 
\begin{equation}
N_+(0)=\chi_+(0)=\int_{M_0}\,{\Td(M_0)\,\Ch(L_0)}=
\RR(M_0,L_0),\end{equation}
q.e.d.

\end{document}